\documentclass[11pt]{article}
\usepackage[tbtags]{amsmath}
\usepackage{amssymb}
\usepackage{amsthm}
\usepackage[misc]{ifsym}
\usepackage{cases}
\usepackage{mathrsfs}
\usepackage{graphicx}
\usepackage{color}
\numberwithin{equation}{section}   
\normalsize
\title{\bf Indefinite linear-quadratic optimal control of mean-field stochastic differential equation with jump diffusion: an  equivalent cost functional method \thanks{This work was supported by the National Natural Science Foundation of China under Grants 61925306, 61821004, and 11831010.}}
\author{\normalsize Guangchen Wang\thanks{\it School of Control Science and Engineering, Shandong University, Jinan 250061, P.R. China,  E-mail: wguangchen@sdu.edu.cn}, Wencan Wang\thanks{\it Corresponding author, School of Control Science and Engineering, Shandong University, Jinan 250100, P.R. China, E-mail: wwencan@163.com}}

\newtheorem{thm}{Theorem}[section]
\newtheorem{definition}{Definition}[section]
\newtheorem{lem}{Lemma}[section]
\newtheorem{rmk}{Remark}[section]
\begin{document}
\maketitle

\noindent{\bf Abstract:}\quad In this paper, we consider a linear-quadratic  optimal control problem of mean-field stochastic differential equation with jump diffusion, which is also called as  an  MF-LQJ problem. Here,   cost functional is allowed to be indefinite. We use an equivalent cost functional method to deal with the MF-LQJ problem with indefinite weighting matrices. Some equivalent cost functionals enable us to establish a bridge between indefinite and positive-definite   MF-LQJ problems. With such a bridge,  solvabilities of stochastic Hamiltonian system  and Riccati equations are further characterized.  Optimal control of the indefinite MF-LQJ problem is represented as a state feedback via solutions of Riccati equations. As a by-product, the method  provides a new way to prove the existence and uniqueness of solution to mean field forward-backward stochastic differential equation with jump diffusion (MF-FBSDEJ, for short), where    existing methods in literature   do not work. Some examples are provided to illustrate our results.

\vspace{2mm}


\vspace{2mm}


\section{Introduction}

In  recent years, there is an increasing interest in mean-field control theory in  mathematics, engineering and finance. Comparing with   classical stochastic  optimal control, a new feature of this problem is that both  objective functional and  system dynamics involve   states and   controls as well as their expected values. There are rich literatures  on deriving  necessary conditions for  optimality. See, for example, Andersson and Djehiche \cite{Daniel2011A}, Buckdahn et al. \cite{BuckdahnA}, Djehiche et al. \cite{djehiche2015stochastic},    Wang et al.  \cite{WangStochastic}. Linear quadratic (LQ, for short) problems of  mean-field type have also been  investigated. Yong \cite{yong2013linear} systematically studied  an LQ problem of  mean-field stochastic differential equation (MF-SDE). Later on,  Sun \cite{sun2017mean} and Li et al. \cite{li2016mean} concerned   open-loop   and closed-loop solvabilities  for LQ problems of MF-SDE, respectively. Elliott et  al. \cite{elliott2013discrete} 
dealt with   an LQ problem of MF-SDE with discrete time  setting. Qi et al. \cite{qi2019stabilization} investigated   stabilization and control
problems for linear MF-SDE  under
standard assumptions. Barreiro-Gomez et al. \cite{barreiro2020discrete, barreiro2019linear} considered LQ mean-field-type games.

It is well known that jump diffusion processes play an increasing role in describing  stochastic dynamical systems, due to its wide application  in  financial, economic and engineering problems. For example, a geometric Brownian motion is usually used to model stock price, but it cannot reflect    discontinuous characteristics, which may be induced by large fluctuations. There are various literatures  on related topic.  The interested readers may refer to  Haadem et al. \cite{Haadem2013infinite},
\O ksendal and Sulem \cite{oksendal2007applied}, Shen et al. \cite{shen2014jump}    for more information. Haadem et al. \cite{Haadem2013infinite} obtained a maximum principle for jump diffusion process with infinite horizon and dealt an optimal portfolio selection problem.
Shen et al. \cite{shen2014jump} investigated stochastic maximum principle of mean field jump diffusion process with delay and applied their results to a bicriteria mean-variance portfolio selection problem. \O ksendal and Sulem \cite{oksendal2007applied} systematically discussed   optimal control, optimal stopping,
and impulse control of   jump diffusion processes.

In  this paper, a kind of indefinite LQ problem of mean field type jump diffusion process   is investigated. Indefinite LQ problems, first studied by Chen et  al. \cite{chen1998stochastic},  have received considerable attention.   Chen et al. \cite{Chen2000Stochastic}  employed the method of completion of squares to study an  indefinite stochastic LQ problem. Huang and Yu \cite{huang2014solvability} and Yu \cite{yu2013equivalent} proposed an equivalent cost functional method to deal with stochastic LQ problem with indefinite weighting matrices.
Ni et al. \cite{ni2015indefinite,ni2014indefinite} considered  indefinite LQ problems of discrete-time MF-SDE for an infinite horizon and
a finite horizon, respectively. Li et al. \cite{2020indefinite} studied an indefinite LQ problem of MF-SDE by introducing a relax  compensator. Wang et al. \cite{wang2020mean} concerned with uniform stabilization and social optimality for general mean field LQ control systems, where state weight is not assumed with the definiteness condition. Indefinite MF-LQJ   problems, which are   natural generalizations of those in  \cite{TangLinear} and \cite{tang2020jump}, have been not yet completely studied. Tang and Meng  \cite{TangLinear} investigated a definite MF-LQJ problem in finite horizon and derived two  Riccati equations by decoupling   optimality system. It was shown that under Assumption (S) given in Section 3 with  $S(t)\equiv 0, \bar S(t)\equiv 0$, these two Riccati equations are uniquely solvable and a feedback representation for optimal control is obtained. We point out that Assumption (S) is exactly the definite condition when we study an MF-LQJ problem. The   MF-LQJ  problem reduces to an  NC-LQJ  problem if $S(t)\equiv 0, \bar S(t)\equiv 0$, where ``NC" is the capital initials for ``no cross". Tang et al. \cite{tang2020jump} focused on    open-loop   and closed-loop solvabilities  of an MF-LQJ problem, which extended results in \cite{sun2017mean,li2016mean}. However, the solvabilities of related Riccati equations without Assumption (S) have not been specified.

Inspired by  \cite{huang2014solvability} and  \cite{yu2013equivalent}, we use an  equivalent cost functional method  to deal with an  indefinite MF-LQJ  problem. As a preliminary result, we discuss a definite MF-LQJ problem. As mentioned above, the results obtained in  \cite{TangLinear} are not applicable for solving this MF-LQJ problem. We introduce an invertible linear transformation, which links the  MF-LQJ problem with the corresponding NC-LQJ problem. Combining the results in existing literature with this linear transformation, we obtain an optimal control of the  MF-LQJ problem under Assumption (S).
Then we introduce two auxiliary functions to construct equivalent cost functionals. The original MF-LQJ problem with indefinite control weighting matrices is transformed into an  MF-LQJ problem under Assumption (S). In a word, we can investigate the indefinite MF-LQJ problem by using this method.

Our paper distinguishes from existing literature in the following aspects.  (i) An indefinite MF-LQJ problem is discussed in this paper, which generalized the results in \cite{TangLinear,tang2020jump,2020indefinite}. The considered model could characterize
more general problems and the jump diffusion item   is important in some controlled dynamics system. As we will see in Example 5.1, there is no  equivalent cost functional satisfying Assumption (S) if jump diffusion item  disappears, which implies that we can not construct an optimal control directly in terms of   Riccati
equations. We further discussed  the existence and uniqueness of solutions to the corresponding stochastic Hamiltonian system and Riccati equations without Assumption (S), which have not been considered in \cite{TangLinear,tang2020jump}. (ii) Compared with   derivation of optimal control for the  definite  NC-LQJ problem in \cite{TangLinear},  we derive an feedback control of ``Problem MF'' under Assumption (S) through a simple calculation. Actually, we introduce  an invertible linear transformation, which links MF-LQJ problem with the corresponding NC-LQJ problem. (iii) Our results  provide  an alternative and effective way to obtain the solvability of an  MF-FBSDEJ, which does not satisfy   classical conditions in existing literature. In fact, when we consider two equivalent cost functionals with the same control system, we can get the equivalence by an invertible linear transformation  between the corresponding stochastic Hamiltonian systems. We point out that the equivalence    is existed in a family of stochastic Hamiltonian systems. 
Therefore, we can prove the solvability of a  more general MF-FBSDEJ. Moreover,  sometimes an MF-FBSDEJ  may coincide  with  the stochastic Hamiltonian system of an MF-LQJ problem, which implies that  the solvability of MF-FBSDEJ is actually the solvability of   corresponding   stochastic Hamiltonian system.  Thus, in order to obtain the unique  solvability of an MF-FBSDEJ, we need only to find an equivalent cost functional satisfying  Assumption (S) of the related MF-LQJ problem. (iv) Relying on the equivalent cost functional method, we can investigate the solvability of indefinite Riccati equations by virtue of the solvability of positive definite Riccati equations. In fact, the original MF-LQJ problem with indefinite control weighting matrices can be transformed into a definite MF-LQJ problem by looking for a simpler and  more flexible equivalent cost functional. And there exists an equivalent relation between the corresponding Riccati equations. Similarly,  we can get   solvabilities of Riccati equations with indefinite condition.

The rest of this paper is organized as follows. In Section 2, we formulate an MF-LQJ
problem and give some assumptions throughout this paper. Section 3 aims to study the
MF-LQJ problem under Assumption (S). We reduce a general MF-LQJ problem to
an NC-LQJ problem via an invertible linear transformation. In Section 4, we present
our main results. We use the equivalent cost functional method to study an
MF-LQJ problem with indefinite weighting matrices. Section 5 gives several illustrative examples. Finally, in Section 6, we conclude this
paper.

\section{Problem formulation}
Let $\mathbb R^{n\times m}$ be an Euclidean space of all $n\times m$ real matrices with inner product $\langle\cdot,\cdot\rangle$ being given by $\langle M,N\rangle\mapsto tr(M^\top N)$, where the superscript $\top$ denotes the transpose of vectors or matrices. The induced norm is given by $|M|=\sqrt{tr(M^\top M)}$. In particular, we denote by $\mathbb S^n$ the set  of all $n\times n$ symmetric matrices. We mean by an $n\times n$ matrix $N\geq 0$ that $N$ is a nonnegative matrix. Let $T>0$  be a fixed time horizon and $(\Omega,\mathcal F,\mathbb F,\mathbb P)$ be a complete filtered probability space. The filtration $\mathbb F\equiv\{\mathcal F_t\}_{t\geq 0}$  is   generated by the following two mutually independent processes,  augmented by all the $\mathbb{P}$-null sets:
a standard 1-dimensional Brownian motion $W_t$ and a Poisson random measure $N(dt, d\theta)$ on $\mathbb R_+\times\Theta$, where $\Theta\subseteq\mathbb R\setminus\{0\}$  is a nonempty  set,   with   compensator   $\widehat{N}(dt, d\theta)=\nu(d\theta)dt$, such that $\widetilde N([0,t], A)=N([0,t], A)-\widehat{N}([0,t], A)$ is a   martingale  for all $A\in  \mathcal{B}(\Theta)$ satisfying  $\nu(A)<\infty$.  $\mathcal{B}(\Theta)$    is the   Borel $\sigma$-field generated by  $\Theta$. 
Here,
$\nu(d\theta)$ is a   $\sigma$-finite measure on $(\Theta, \mathcal{B}(\Theta))$ satisfying $\int_{\Theta}(1\wedge \theta^2)\nu(d\theta)<\infty$, which is called the characteristic measure. Then    $\widetilde N(dt, d\theta)=N(dt, d\theta)-\nu(d\theta)dt$ is the compensated Poisson random measure. 
For any Euclidean space $M$, we introduce the following spaces:\\
$L^\infty(0,T;M)=\Big\{u:[0,T]\to M|u(\cdot)$ is  a   bounded function\Big\};\\
$L_{\mathbb F}^2(0,T;M)=\Big\{u:[0,T]\times \Omega\to M|u(\cdot)$ is an $\mathbb F$-adapted stochastic  process such that $\mathbb E\left[\int_0^T|u(t)|^2dt\right]<\infty\Big\};$\\
$S_{\mathbb F}^2(0,T;M)=\Big\{u:[0,T]\times \Omega\to  M| u(\cdot)$ is  an  $\mathbb F$-adapted $c\grave{a}dl\grave{a}g$ process in $L_{\mathbb F}^2(0,T;M)$ such that $\mathbb E\big[\sup_{t\in[0,T]}|u(t)|^2\big]<\infty\Big\}$;\\
$L_{\nu}^2(M)=\Big\{r:[0,T]\times \Theta\to M|r(\cdot)$ is  a   deterministic  function such that $sup_{t\in [0,T]}\int_{\Theta}| r(t,\theta)|^2\nu(d\theta)<\infty\Big\}$;\\
$L_{\mathbb F,\nu}^2(0,T;M)=\bigg\{r:[0,T]\times \Theta\times \Omega\to M| r(\cdot)$ is an $\mathbb F$-predictable stochastic process such that  $\mathbb E\left[\int_0^T\int_{\Theta}|r(t,\theta)|^2\nu(d\theta)dt\right] <\infty\bigg\}$.\\

 Consider a controlled linear MF-SDEJ
 \begin{equation}\label{station 1}
\left\{ \begin{aligned}
dX_{t}=&\Big\{A_{t}X_{t}+\bar A_{t}\mathbb E[X_{t}]+B_{t}u_{t}+\bar B_{t}\mathbb E[u_{t}]\Big\}dt\\&+\Big\{C_{t}X_{t}+\bar C_{t}\mathbb E[X_{t}]+D_{t}u_{t}+\bar D_{t}\mathbb E[u_{t}]\Big\}dW_{t}\\&+\int_{\Theta}\Big\{E_{t,\theta}X_{t-}+\bar E_{t,\theta}\mathbb E[X_{t-}]+F_{t,\theta}u_{t}+\bar F_{t,\theta}\mathbb E[u_{t}]\Big \}\widetilde N(dt,d\theta),\qquad t\in{[0,T]},\\
X_{0}=&\ x\in \mathbb R^n,
\end{aligned}\right.
\end{equation}
where $A_{t}, \bar A_{t},B_{t},\bar B_{t},C_{t},\bar C_{t},D_{t},\bar D_{t}, E_{t,\theta}, \bar E_{t,\theta},F_{t,\theta}, \bar F_{t,\theta}$ are given matrix valued deterministic functions. In the above equation, $u$, valued in $\mathbb R^m$, is a control process and $X$, valued in $\mathbb R^n$, is the corresponding state process.
In this paper, an  admissible control $u$ is defined as a predictable process such that $u\in L_{\mathbb F}^2(0,T;\mathbb R^m)$. The set of all   admissible controls   is denoted by  $\mathcal U[0,T]$.
We introduce a cost functional
\begin{equation}\label{cost functional}
\begin{aligned}
J[u]=&\ \frac{1}{2}\mathbb E\left\{\langle GX_{T},X_{T}\rangle+\langle \bar G\mathbb E[X_{T}],\mathbb E[X_{T}]\rangle\right.\\&+\int_0^T\left\langle\left(\begin{array}{cc}
Q_{t}&S_{t}\\S_{t}^\top&R_{t}
\end{array}\right)\left(\begin{array}{c}
X_{t}\\u_{t}
\end{array}\right),\left(\begin{array}{c}
X_{t}\\u_{t}
\end{array}\right)\right\rangle dt\\&+\left.\int_0^T\left\langle\left(\begin{array}{cc}
\bar Q_{t}&\bar S_{t}\\ \bar S_{t}^\top&\bar R_{t}
\end{array}\right)\left(\begin{array}{c}
\mathbb E[X_{t}]\\\mathbb E[u_{t}]
\end{array}\right),\left(\begin{array}{c}
\mathbb E[X_{t}]\\\mathbb E[u_{t}]
\end{array}\right)\right\rangle dt\right\},
\end{aligned}
\end{equation}
where $G,\bar G$ are symmetric matrices and $Q_{t}, S_{t}, R_{t}, \bar Q_{t}, \bar S_{t},\bar R_{t}$ are deterministic matrix-valued functions with $Q_{t}=Q_{t}^\top, R_{t}=R_{t}^\top,\bar Q_{t}=\bar Q_{t}^\top,\bar R_{t}=\bar R_{t}^\top$. Our MF-LQJ   problem is stated as follows.\\
{\it\textbf{Problem MF}}: For any $x\in \mathbb R^n$,  find a $u^*\in \mathcal U[0,T]$ such that
\begin{equation}\label{inf}
J[u^*]=\inf_{u\in \mathcal U[0,T]}J[u].
\end{equation}
Any $u^*\in \mathcal U[0,T]$ satisfying \eqref{inf} is called an optimal control of Problem MF, and the corresponding state process $X^*=X(x,u^*)$ is called an optimal state process. $(X^*, u^*)$ is called an optimal pair.

The following assumptions will be in force throughout this paper.
\\(H1) The coefficients of  state equation satisfy
\begin{equation*}
\left\{\begin{array}{l}
A, \bar A, C, \bar C\in L^\infty(0,T;\mathbb R^{n\times n}),
B, \bar B, D, \bar D\in L^\infty(0,T;\mathbb R^{n\times m}),\\
E, \bar E\in L_{\nu}^2(\mathbb R^{n\times n}), F, \bar F\in L_{\nu}^2(\mathbb R^{n\times m}).
\end{array}\right.
\end{equation*}
(H2) The weighting matrices in    cost functional satisfy
\begin{equation*}
\left\{\begin{aligned}
&Q, \bar Q \in L^\infty(0,T;\mathbb S^{n}),\  R, \bar R \in L^\infty(0,T;\mathbb S^{m}),\\
&S, \bar S\in L^\infty(0,T;\mathbb R^{n\times m}), \ G, \bar G\in \mathbb S^{n}.
\end{aligned}\right.
\end{equation*}
We can show that under (H1), for any $u\in \mathcal U[0,T]$, \eqref{station 1} admits a unique solution $X=X(x,u)\in S_{\mathbb F}^2(0,T;\mathbb R^n)$.

\section{MF-LQJ problem under standard conditions}
In this section, we aim at studying Problem MF under some standard conditions. We introduce   an invertible linear transformation, which links MF-LQJ problem with the corresponding NC-LQJ problem.\\
{\it\textbf{Assumption (S)}}:  For some $\alpha_0>0$,
\begin{equation*}
\left\{\begin{aligned}
&R_{t}, R_{t}+\bar R_{t}\geq \alpha_0 I, Q_{t}-S_{t}R_{t}^{-1}S_{t}^\top\geq 0, \\&Q_{t}+\bar Q_{t}-(S_{t}+\bar S_{t})(R_{t}+\bar R_{t})^{-1}(S_{t}+\bar S_{t})^\top\geq 0, \quad t\in[0,T],
\\&G, G+\bar G\geq 0.
\end{aligned}
\right.
\end{equation*}

We introduce  a stochastic Hamiltonian system related to Problem MF
\begin{equation}\label{H   system}
\left\{\begin{aligned}&dX_{t}=\Big\{A_{t}X_{t}+\bar A_{t}\mathbb E[X_{t}]+B_{t}u_{t}+\bar B_{t}\mathbb E[u_{t}]\Big\}dt
\\&\qquad\ \ \ +\Big\{C_{t}X_{t}+\bar C_{t}\mathbb E[X_{t}]+D_{t}u_{t}+\bar D_{t}\mathbb E[u_{t}]\Big\}dW_t\\&\qquad\ \ \ +\int_{\Theta}\Big\{E_{t,\theta}X_{t-}+\bar E_{t,\theta}\mathbb E[X_{t-}]+F_{t,\theta}u_{t}+\bar F_{t,\theta}\mathbb E[u_{t}]\Big \}\widetilde N(dt,d\theta),\\&dY_{t}=-\Big\{A_{t}^\top Y_{t}+\bar A_{t}^\top \mathbb E[Y_{t}]+C_{t}^\top Z_{t}+\bar C_{t}^\top\mathbb E[Z_{t}]+\int_{\Theta} \Big(E_{t,\theta}^\top r_{t,\theta}+\bar E_{t,\theta}^\top \mathbb E[r_{t,\theta}]\Big)\nu(d\theta)\\&\qquad \ \ \ +Q_{t}X_{t}+\bar Q_{t}\mathbb E[X_{t}]+S_{t}u_{t}+\bar S_{t}\mathbb E[u_{t}]\Big\}dt+Z_{t}dW_{t}+\int_{\Theta} r_{t,\theta}\widetilde N(dt,d\theta),\\&X_0=x,\ \ \ \  Y_{T}=GX_{T}+\bar G\mathbb E[X_{T}],\\&R_{t}u_{t}+\bar R_{t}\mathbb E[u_{t}]+S_{t}^\top X_{t-}+\bar S_{t}^\top\mathbb E[X_{t-}]+B_{t}^\top Y_{t-}+\bar B_{t}^\top\mathbb E[Y_{t-}]\\&\ \ +D_{t}^\top Z_{t}+\bar D_{t}^\top\mathbb E[Z_{t}]+\int_{\Theta} \left(F_{t,\theta}^\top r_{t,\theta}+\bar F_{t,\theta}^\top \mathbb E[r_{t,\theta}]\right)\nu(d\theta)=0.
\end{aligned}
\right.
\end{equation}
Using the method in \cite{TangLinear}, we decouple the above Hamiltonian  system and derive  Riccati equations associated with Problem MF
\begin{equation}\label{P riccati}
\left\{\begin{aligned}
&\dot{P_{t}}+P_{t}A_{t}+A_{t}^\top P_{t}+C_{t}^\top P_{t}C_{t}+\int_{\Theta} E_{t,\theta}^\top P_{t}E_{t,\theta}\nu(d\theta) +Q_{t}
\\&\  -\left(S_{t}+P_{t}B_{t}+C_{t}^\top P_{t}D_{t}+\int_{\Theta} E_{t,\theta}^\top P_{t} F_{t,\theta}\nu(d\theta)\right)\Sigma_{0t}^{-1}
\\&\ \cdot \left(S_{t}^\top +B_{t}^\top P_{t}+D_{t}^\top P_{t}C_{t}+\int_{\Theta} F_{t,\theta}^\top P_{t} E_{t,\theta}\nu(d\theta)\right)=0,\\&P_{T}=G,
\end{aligned}\right.
\end{equation}
\begin{equation}\label{Pi riccati}
\left\{\begin{aligned}
&\dot{\Pi_{t}}+\Pi_{t}(A_{t}+\bar A_{t})+(A_{t}+\bar A_{t})^\top \Pi_{t}+(C_{t}+\bar C_{t})^\top P_{t}(C_{t}+\bar C_{t})
\\&\ +\int_{\Theta} (E_{t,\theta}+\bar E_{t,\theta})^\top P_{t}(E_{t,\theta}+\bar E_{t,\theta})\nu(d\theta)+Q_{t}+\bar Q_{t}
\\&\ -\left[(S_{t}+\bar S_{t})+\Pi_{t}(B_{t}+\bar B_{t})+(C_{t}+\bar C_{t})^\top P_{t}(D_{t}+\bar D_{t})\right.
\\&\ \left.+\int_{\Theta} (E_{t,\theta}+\bar E_{t,\theta})^\top P_{t} (F_{t,\theta}+\bar F_{t,\theta})\nu(d\theta)\right]\Sigma_{1t}^{-1}
\\&\ \cdot\left[(S_{t}+\bar S_{t})^\top+(B_{t}+\bar B_{t})^\top \Pi_{t}+(D_{t}+\bar D_{t})^\top P_{t}(C_{t}+\bar C_{t})\right.
\\&\ +\left.\int_{\Theta} (F_{t,\theta}+F_{t,\theta})^\top P_{t} (E_{t,\theta}+E_{t,\theta})\nu(d\theta)\right]= 0,\\ &\Pi_{T}=G+\bar G,
\end{aligned}
\right.
\end{equation}
where
\begin{equation*}
\begin{aligned}
&\Sigma_{0t}=R_{t}+D_{t}^\top P_{t}D_{t}+\int_{\Theta} F_{t,\theta}^\top P_{t} F_{t,\theta}\nu(d\theta),
\\&\Sigma_{1t}=R_{t}+\bar R_{t}+(D_{t}+\bar D_{t})^\top P_{t}(D_{t}+\bar D_{t})+\int_{\Theta} (F_{t,\theta}+F_{t,\theta})^\top P_{t} (F_{t,\theta}+F_{t,\theta})\nu(d\theta).
\end{aligned}
\end{equation*}
Note that  \eqref{H   system} is a coupled MF-FBSDEJ, where the coupling comes from the last relation (which is essentially the maximum condition in the usual Pontryagin type maximum principle).
Different from an  NC-LQJ  problem, there are   additional items $2\langle X_{t},S_{t}u_{t}\rangle$ and  $2\langle\mathbb E[X_{t}],\bar S_{t}\mathbb E[u_{t}]\rangle$   in  cost functional \eqref{cost functional}. Next, we want to  reduce Problem MF to an NC-LQJ problem. For this, we introduce a controlled system
\begin{equation}\label{NC state}
\left\{ \begin{aligned}
d\tilde X_{t}=
&\ \Big\{A_{1t}(\tilde X_{t}-\mathbb E[\tilde X_{t}])+B_{t}\tilde u_{t}+[A_{1t}+\bar A_{1t}]\mathbb E[\tilde X_{t}]+\bar B_{t}\mathbb E[\tilde u_{t}]\Big\}dt
\\&+\Big\{C_{1t}(\tilde X_{t}-\mathbb E[\tilde X_{t}])+D_{t}\tilde u_{t}+[C_{1t}+\bar C_{1t}]\mathbb E[\tilde X_{t}]+\bar D_{t}\mathbb E[\tilde u_{t}]\Big\}dW_{t}
\\&+\int_{\Theta}\Big\{E_{1t,\theta}(\tilde X_{t-}-\mathbb E[\tilde X_{t-}])+F_{t,\theta}\tilde u_{t}+[E_{1t,\theta}+\bar E_{1t,\theta}]\mathbb E[\tilde X_{t-}]+\bar F_{t,\theta}\mathbb E[\tilde u_{t}]\Big \}\widetilde N(dt,d\theta),\\
\tilde X_{0}=\ & x,
\end{aligned}\right.
\end{equation}
and a cost functional
\begin{equation}\label{NC cost}
\begin{aligned}
&\tilde J[\tilde u]=\frac{1}{2}\mathbb E\left\{\langle G\tilde X_{T},\tilde X_{T}\rangle+\langle \bar G\mathbb E[\tilde X_{T}],\mathbb E[\tilde X_{T}]\rangle\right.\\&\ \ +\int_0^T\left\langle\left(\begin{array}{cc}
Q_{1t}&\mathbf{0}_{n\times m}\\ \mathbf{0}_{m\times n}&R_{t}
\end{array}\right)\left(\begin{array}{c}
\tilde X_{t}-\mathbb{E}[\tilde X_{t}]\\ \tilde u_{t}-\mathbb{E}[\tilde u_{t}]
\end{array}\right),\left(\begin{array}{c}
\tilde X_{t}-\mathbb{E}[\tilde X_{t}]\\ \tilde u_{t}-\mathbb{E}[\tilde u_{t}]
\end{array}\right)\right\rangle dt
\\&\ \ \left.+\int_0^T\left\langle\left(\begin{array}{cc}
Q_{1t}+\bar Q_{1t}&\mathbf{0}_{n\times m}\\\mathbf{0}_{m\times n}&R_{t}+\bar R_{t}
\end{array}\right)\left(\begin{array}{c}
\mathbb E[\tilde X_{t}]\\\mathbb E[\tilde u_{t}]
\end{array}\right),\left(\begin{array}{c}
\mathbb E[\tilde X_{t}]\\\mathbb E[\tilde u_{t}]
\end{array}\right)\right\rangle dt\right\},
\end{aligned}
\end{equation}
where
\begin{equation*}
\begin{aligned}
&A_{1t}=A_{t}-B_{t}R_{t}^{-1}S_{t}^\top,
\\&A_{1t}+\bar A_{1t}=A_{t}+\bar A_{t}-(B_{t}+\bar B_{t})(R_{t}+\bar R_{t})^{-1}(S_{t}+\bar S_{t})^\top,
\\
&C_{1t}=C_{t}-D_{t}R_{t}^{-1}S_{t}^\top,
\\&C_{1t}+\bar C_{1t}=C_{t}+\bar C_{t}-(D_{t}+\bar D_{t})(R_{t}+\bar R_{t})^{-1}(S_{t}+\bar S_{t})^\top,
\\
&E_{1t,\theta}=E_{t,\theta}-F_{t,\theta}R_{t}^{-1}S_{t}^\top,
\\&E_{1t,\theta}+\bar E_{1t,\theta}=E_{t,\theta}+\bar E_{t,\theta}-(F_{t,\theta}+\bar F_{t,\theta})(R_{t}+\bar R_{t})^{-1}(S_{t}+\bar S_{t})^\top,
\\
&Q_{1t}=Q_{t}-S_{t}R_{t}^{-1}S_{t}^\top,
\\&Q_{1t}+\bar Q_{1t}=Q_{t}+\bar Q_{t}-(S_{t}+\bar S_{t})(R_{t}+\bar R_{t})^{-1}(S_{t}+\bar S_{t})^\top.
\end{aligned}
\end{equation*}
The corresponding NC-LQJ problem is  stated as follows.

{\it\textbf{Problem NC}}: For any $x\in \mathbb R^n$,  find a $\tilde u^*\in \mathcal U[0,T]$ such that
\begin{equation}\label{NC inf}
\tilde J[\tilde u^*]=\inf_{\tilde u\in \mathcal U[0,T]}\tilde J[\tilde u].
\end{equation}

Similar to Problem MF, we write the stochastic Hamiltonian system and   Riccati equations corresponding to Problem NC
 \begin{equation}\label{H  system NC}
\left\{\begin{aligned}&d\tilde X_{t}=\left\{A_{1t}\tilde X_{t}+\bar A_{1t}\mathbb E[\tilde X_{t}]+B_{t}\tilde u_{t}+\bar B_{t}\mathbb E[\tilde u_{t}]\right\}dt\\&\ \ \ \ \ \ \ +\left\{C_{1t}\tilde X_{t}+\bar C_{1t}\mathbb E[\tilde X_{t}]+D_{t}\tilde u_{t}+\bar D_{t}\mathbb E[\tilde u_{t}]\right\}dW_{t}\\&\ \ \ \ \ \ \ +\int_{\Theta}\Big\{E_{1t,\theta}\tilde X_{t-}+\bar E_{1t,\theta}\mathbb E[\tilde X_{t-}]+F_{t,\theta}\tilde u_{t}+\bar F_{t,\theta}\mathbb E[\tilde u_{t}]\Big \}\widetilde N(dt,d\theta),\\&d\tilde Y_{t}=-\Big\{A_{1t}^\top\tilde Y_{t}+\bar A_{1t}^\top\mathbb E[\tilde Y_{t}]+C_{1t}^\top\tilde Z_{t}+\bar C_{1t}^\top\mathbb E[\tilde Z_{t}]+\int_{\Theta} \Big(E_{1t,\theta}^\top \tilde{r}_{t,\theta}+\bar E_{1t,\theta}^\top \mathbb E[\tilde{r}_{t,\theta}]\Big)\nu(d\theta)\\&\ \ \ \ \ \ \ +Q_{1t}\tilde X_{t}+\bar Q_{1t}\mathbb E[\tilde X_{t}]\Big\}dt +\tilde Z_{t}dW_{t}+\int_{\Theta} \tilde r_{t,\theta}\widetilde N(dt,d\theta),\\&\tilde X_{0}=x, \ \ \ \  \tilde Y_{T}=G\tilde X_{T}+\bar G\mathbb E[\tilde X_{T}],\\&R_{t}\tilde u_{t}+\bar R_{t}\mathbb E[\tilde u_{t}]+B_{t}^\top \tilde Y_{t-}+\bar B_{t}^\top\mathbb E[\tilde Y_{t-}]+D_{t}^\top \tilde Z_{t}+\bar D_{t}^\top\mathbb E[\tilde Z_{t}]\\&+\int_{\Theta} \left(F_{t,\theta}^\top \tilde{r_{t,\theta}}+\bar F_{t,\theta}^\top \mathbb E[\tilde{r_{t,\theta}}]\right)\nu(d\theta)=0,
\end{aligned}
\right.
\end{equation}

\begin{equation}\label{P riccati NC}
\left\{\begin{aligned}
&\dot{P_{t}}+P_{t}A_{1t}+A_{1t}^\top P_{t}+ C_{1t}^\top P_{t} C_{1t}+\int_{\Theta} E_{1t,\theta}^\top P_{t}E_{1t,\theta}\nu(d\theta)+ Q_{1t}
\\&\quad-\left(P_{t}B_{t}+ C_{1t}^\top P_{t}D_{t}+\int_{\Theta} E_{1t,\theta}^\top P_{t} F_{t,\theta}\nu(d\theta)\right)\Sigma_{0t}^{-1}
\\&\quad\cdot\left(B_{t}^\top  P_{t}+D_{t}^\top P_{t}C_{1t}+\int_{\Theta} F_{t,\theta}^\top P_{t} E_{1t,\theta}\nu(d\theta)\right)=0,\\&P_{T}=G,
\end{aligned}\right.
\end{equation}
\begin{equation}\label{Pi riccati NC}
\left\{\begin{aligned}
&\dot{\Pi_{t}}+\Pi_{t}\left(A_{1t}+\bar A_{1t}\right)+\left(A_{1t}+\bar A_{1t}\right)^\top \Pi_{t}+\left(C_{1t}+\bar C_{1t}\right)^\top P_{t}\left(C_{1t}+\bar C_{1t}\right)
\\&\quad+\int_{\Theta} (E_{1t,\theta}+\bar E_{1t,\theta})^\top P_{t}(E_{1t,\theta}+\bar E_{1t,\theta})\nu(d\theta)+Q_{1t}+\bar Q_{1t}
\\&\quad-\left[\Pi_{t}(B_{t}+\bar B_{t})+(C_{1t}+\bar C_{1t})^\top P_{t}(D_{t}+\bar D_{t})+\int_{\Theta} (E_{1t,\theta}+\bar E_{1t,\theta})^\top P_{t}(F_{t,\theta}+\bar F_{t,\theta})\nu(d\theta)\right]\Sigma_{1t}^{-1}
\\&\quad\cdot\left[(B_{t}+\bar B_{t})^\top \Pi_{t}+(D_{t}+\bar D_{t})^\top P_{t}(C_{1t}+\bar C_{1t})+\int_{\Theta} (F_{t,\theta}+F_{t,\theta})^\top P_{t}(E_{1t,\theta}+E_{1t,\theta})\nu(d\theta)\right]= 0, \\ &\Pi_{T}=G+\bar G,
\end{aligned}
\right.
\end{equation}
where
\begin{equation*}
\begin{aligned}
&\Sigma_{0t}=R_{t}+D_{t}^\top P_{t}D_{t}+\int_{\Theta} F_{t,\theta}^\top P_{t} F_{t,\theta}\nu(d\theta),
\\&\Sigma_{1t}=R_{t}+\bar R_{t}+(D_{t}+\bar D_{t})^\top P_{t}(D_{t}+\bar D_{t})+\int_{\Theta} (F_{t,\theta}+F_{t,\theta})^\top P_{t} (F_{t,\theta}+F_{t,\theta})\nu(d\theta).
\end{aligned}
\end{equation*}
\begin{lem}\label{control system lemma}
Let   Assumption (S)  hold. For any two pairs $(X,u)$ and $(\tilde X,\tilde u)$, we introduce a linear transformation
\begin{equation}\label{linear   trans1}
\left(\begin{array}{c}
X-\mathbb E [X]\\
u-\mathbb E [u]\end{array}\right)=\left(\begin{array}{cc}
I_{n\times n}&\mathbf{0}_{n\times m}\\
-R^{-1}S^\top&I_{m\times m}\end{array}\right)\left(\begin{array}{c}
\tilde X-\mathbb E [\tilde X]\\
\tilde u-\mathbb E [\tilde u]\end{array}\right),
\end{equation}
 \begin{equation}\label{linear   trans2}
\left(\begin{array}{c}
\mathbb E [X]\\
\mathbb E[u]\end{array}\right)=\left(\begin{array}{cc}
I_{n\times n}&\mathbf{0}_{n\times m}\\
-(R+\bar R)^{-1}(S+\bar S)^\top&I_{m\times m}\end{array}\right)\left(\begin{array}{c}
\mathbb E [\tilde X]\\
\mathbb E [\tilde u]\end{array}\right).
\end{equation}
Then the following two statements are equivalent:
\begin{itemize}
\item [(i).]$(X,u)$ is an admissible (optimal) control of {\it Problem MF}.
\item [(ii).]$(\tilde X,\tilde u)$ is an admissible (optimal) control of {\it Problem NC}.
\end{itemize}
Moreover, we have $J[u]=\tilde J[\tilde u].$
\end{lem}
\begin{proof}
It follows from \eqref{linear   trans1} and \eqref{linear   trans2} that
\begin{equation*}
\left(\begin{array}{c}
\tilde X-\mathbb E [\tilde X]\\
\tilde u-\mathbb E [\tilde u]\end{array}\right)=\left(\begin{array}{cc}
I_{n\times n}&\mathbf{0}_{n\times m}\\
R^{-1}S^\top&I_{m\times m}\end{array}\right)\left(\begin{array}{c}
X-\mathbb E[X]\\
u-\mathbb E [u]\end{array}\right),
\end{equation*}
 \begin{equation*}
\left(\begin{array}{c}
\mathbb E [\tilde X]\\
\mathbb E [\tilde u]\end{array}\right)=\left(\begin{array}{cc}
I_{n\times n}&\mathbf{0}_{n\times m}\\
(R+\bar R)^{-1}(S+\bar S)^\top&I_{m\times m}\end{array}\right)\left(\begin{array}{c}
\mathbb E [X]\\
\mathbb E [u]\end{array}\right).
\end{equation*}
Then linear transformation \eqref{linear   trans1} with \eqref{linear   trans2} is  invertible. Through direct calculation, it is easy to verify statement  (i)   is equivalent to  statement  (ii) , and thus  $J[u]=\tilde J[\tilde u].$
\end{proof}
The above lemma tells us that there exists some equivalent relationship between
Problem MF and Problem NC. We now analyze the   relationship in terms of
stochastic Hamiltonian system and Riccati equations, respectively.
\begin{lem}\label{H system lemma}
Under  Assumption (S), $(\tilde X^*, \tilde u^*, Y, Z, r)$ is the solution of \eqref{H  system NC} if and only if  $(X^*, u^*, Y, Z, r)$  is the solution of \eqref{H   system}. \end{lem}
\begin{proof}
According to Lemma \ref{control system lemma}, it is not difficult to draw the conclusion.

\end{proof}

\begin{lem}\label{Riccati lemma}
Under Assumption (S), we have
\begin{itemize}
\item[1.] Riccati equations \eqref{P riccati NC} and \eqref{P riccati} are the same.
\item[2.] Riccati equations \eqref{Pi riccati NC} and \eqref{Pi riccati} are the same.
\end{itemize}
\end{lem}
\begin{proof}
For simplicity of notations, we denote
\begin{equation*}
\begin{aligned}
&G(A_{t},B_{t},C_{t},D_{t},E_{t,\theta},F_{t,\theta};Q_{t},R_{t},S_{t};P_{t})
\\=&\ P_{t}A_{t}+A_{t}^\top P_{t}+C_{t}^\top P_{t}C_{t}+\int_{\Theta} E_{t,\theta}^\top P_{t}E_{t,\theta}\nu(d\theta)+Q_{t}
\\&-\left(P_{t}B_{t}+C_{t}^\top P_{t}D_{t}+\int_{\Theta} E_{t,\theta}^\top P_{t} F_{t,\theta}\nu(d\theta)+S_{t}\right)\Sigma_{0t}^{-1}
\\&\cdot\left(B_{t}^\top P_{t}+D_{t}^\top P_{t}C_{t}+\int_{\Theta} F_{t,\theta}^\top P_{t}E_{t,\theta}\nu(d\theta)+S_{t}^\top\right).
\end{aligned}
\end{equation*}
It is enough to prove $G(A_{t},B_{t},C_{t},D_{t},E_{t,\theta},F_{t,\theta};Q_{t},R_{t},S_{t};P_{t}) =G(A_{1t},B_{t},C_{1t},D_{t},E_{1t,\theta},F_{t,\theta};Q_{1t},R_{t},\mathbf 0;P_{t})$.
We have
\begin{equation*}
\begin{aligned}
&\Sigma_{0t}^{-1}\left(B_{t}^\top P_{t}+D_{t}^\top P_{t}C_{t}+\int_{\Theta} F_{t,\theta}^\top P_{t}E_{t,\theta}\nu(d\theta)+S_{t}^\top\right)
\\&-\Sigma_{0t}^{-1}\left(B_{t}^\top P_{t}+D_{t}^\top P_{t}C_{1t}+\int_{\Theta} F_{t,\theta}^\top P_{t}E_{1t,\theta}\nu(d\theta)\right)
\\=&\ \Sigma_{0t}^{-1}\left(D_{t}^\top P_{t}D_{t}+\int_{\Theta} F_{t,\theta}^\top P_{t} F_{t,\theta}\nu(d\theta)+R_{t}\right)R_{t}^{-1}S_{t}^\top
\\=&\ R_{t}^{-1}S_{t}^\top.
\end{aligned}
\end{equation*}

Thus, we calculate
\begin{equation*}
\begin{aligned}
&G(A_{t},B_{t},C_{t},D_{t},E_{t,\theta},F_{t,\theta};Q_{t},R_{t},S_{t};P_{t}) -G(A_{1t},B_{t},C_{1t},D_{t},E_{1t,\theta},F_{t,\theta};Q_{1t},R_{t},\mathbf 0;P_{t})
\\=&\ P_{t}B_{t}R_{t}^{-1}S_{t}^\top+S_{t}R_{t}^{-1}B_{t}^\top P_{t}+C_{t}^\top P_{t}D_{t}R_{t}^{-1}S_{t}^\top+S_{t}R_{t}^{-1}D_{t}^\top P_{t}C_{1t}
\\&+\int_{\Theta} E_{t,\theta}^\top P_{t}F_{t,\theta}R_{t}^{-1}S_{t}^\top\nu(d\theta) +\int_{\Theta} S_{t}R_{t}^{-1}F_{t,\theta}^\top P_{t,\theta}E_{1t,\theta}\nu(d\theta)
\\&+S_{t}R_{t}^{-1}S_{t}^\top-\left(S_{t}+P_{t}B_{t}+C_{t}^\top P_{t}D_{t}+\int_{\Theta}E_{t,\theta}^\top P_{t}F_{t,\theta}\nu(d\theta)\right)R_{t}^{-1}S_{t}^\top
\\&-S_{t}R_{t}^{-1}\left(B_{t}^\top P_{t}+D_{t}^\top P_{t}C_{1t}+\int_{\Theta} F_{t,\theta}^\top P_{t} E_{1t,\theta}\nu(d\theta)\right)
\\=&\ \mathbf{0}_{n\times  n}.
\end{aligned}
\end{equation*}
Consequently, we complete the proof.
\end{proof}

We cite a   result  in  \cite{TangLinear}, which plays a  role in deriving an optimal control of Problem MF under standard condition.
\begin{lem}\label{Yong}
Let (H1)-(H2) and Assumption (S) hold, then  Riccati equations \eqref{P riccati NC} and \eqref{Pi riccati NC} admit unique solutions $P$ and $\Pi$, respectively. Further, the
optimal pair $(\tilde X^*, \tilde u^*)$ of {\it Problem NC} satisfies
\begin{equation*}
\left\{\begin{aligned}
\tilde u_{t}^*=&-\Sigma_{0t}^{-1}\left(B_{t}^\top P_{t}+D_{t}^\top P_{t}C_{1t}+\int_{\Theta} F_{t,\theta}^\top P_{t}E_{1t,\theta}\nu(d\theta)\right)(\tilde X_{t}^*-\mathbb E[\tilde X_{t}^*])\\&-\Sigma_{1t}^{-1}\left[(B_{t}+\bar B_{t})^\top \Pi_{t}+(D_{t}+\bar D_{t})^\top P_{t}(C_{1t}+\bar C_{1t})\right.\\&+\left.\int_{\Theta}(F_{t,\theta}+\bar F_{t,\theta})^\top P_{t}(E_{1t,\theta}+\bar E_{1t,\theta})\nu(d\theta)\right]\mathbb E[\tilde X_{t}^*],\\
d\tilde X_{t}^*=&\ \Big\{A_{1t}\tilde X_{t}^*+\bar A_{1t}\mathbb E[\tilde X_{t}^*]+B_{t}\tilde u_{t}^*+\bar B_{t}\mathbb E[\tilde u_{t}^*]\Big\}dt\\&+\Big\{C_{1t}\tilde X_{t}^*+\bar C_{1t}\mathbb E[\tilde X_{t}^*]+D_{t}\tilde u_{t}^*+\bar D_{t}\mathbb E[\tilde u_{t}^*]\Big\}dW_{t}\\&+\int_{\Theta}\Big\{E_{1t,\theta}\tilde X_{t-}^*+\bar E_{1t,\theta}\mathbb E[\tilde X_{t-}^*]+F_{t,\theta}\tilde u_{t}^*+\bar F_{t,\theta}\mathbb E[\tilde u_{t}^*]\Big \}\widetilde N(dt,d\theta),\\
\tilde X_{0}^*=&\ x.
\end{aligned}
\right.
\end{equation*}
Defining
\begin{equation*}
\left\{\begin{aligned}
\tilde Y_{t}^*=&\ P_{t}(\tilde X_{t}^*-\mathbb E[\tilde X_{t}^*])+\Pi_{t}\mathbb E[\tilde X_{t}^*],
\\ \tilde Z_{t}^*=&\ \left[P_{t}C_{1t}-P_{t}D_{t}\Sigma_{0t}^{-1}\left(B_{t}^\top P_{t}+D_{t}^\top P_{t}C_{1t}+\int_{\Theta} F_{t,\theta}^\top P_{t}E_{1t,\theta}\nu(d\theta)\right)\right]
(\tilde X_{t}^*-\mathbb E[\tilde X_{t}^*])
\\&+\left\{P_{t}(C_{1t}+\bar C_{1t})-P_{t}(D_{t}+\bar D_{t})\Sigma_{1t}^{-1}\left[(B_{t}+\bar B_{t})^\top \Pi_{t}+(D_{t}+\bar D_{t})^\top P_{t}(C_{1t}+\bar C_{1t})\right.\right.
\\&\left.\left.+\int_{\Theta} (F_{t,\theta}+\bar F_{t,\theta})^\top P_{t}(E_{1t,\theta}+\bar E_{1t,\theta})\nu(d\theta)\right]\right\}\mathbb E[\tilde X_{t}^*],\\
\tilde r_{t,\theta}^*=&\ \left[P_{t}E_{1t,\theta}-P_{t}F_{t,\theta}\Sigma_{0t}^{-1}\left(B_{t}^\top P_{t}+D_{t}^\top P_{t}C_{1t}+\int_{\Theta} F_{t,\theta}^\top P_{t}E_{1t,\theta}\nu(d\theta)\right)\right](\tilde X_{t}^*-\mathbb E[\tilde X_{t}^*])
\\&+\Big\{P_{t}(E_{1t,\theta}+\bar E_{1t,\theta})-P_{t}(F_{t,\theta}+\bar F_{t,\theta})\Sigma_{1t}^{-1}\Big[(B_{t}+\bar B_{t})^\top \Pi_{t}+(D_{t}+\bar D_{t})^\top P_{t}(C_{1t}+\bar C_{1t})
\\&\left.\left.+\int_{\Theta} (F_{t,\theta}+\bar F_{t,\theta})^\top P_{t}(E_{1t,\theta}+\bar E_{1t,\theta})\nu(d\theta)\right]\right\}\mathbb E[\tilde X_{t}^*],
\end{aligned}
\right.
\end{equation*}
the $5$-tuple $(\tilde X^*,\tilde u^*, \tilde Y^*,\tilde Z^*,\tilde r^*)$ is the unique   solution to MF-FBSDEJ \eqref{H  system NC}. Moreover,
\begin{equation*}
\inf_{\tilde u \in \mathcal U[0,T]}\tilde J[\tilde u]=\frac{1}{2}\langle\Pi_0x,x\rangle,\ \ \ \forall x\in \mathbb R^n.
\end{equation*}
\end{lem}

Using the   above lemmas, we obtain a main result of this section.
\begin{thm} \label{MF STN}If (H1)-(H2) and Assumption (S) hold, then Riccati equations \eqref{P riccati} and \eqref{Pi riccati} admit unique solutions $P$ and $\Pi$, respectively. Further, the
optimal pair $(X^*, u^*)$ of {\it Problem MF} satisfies
\begin{equation*}
\left\{\begin{aligned}
 u_{t}^*=&-\Sigma_{0t}^{-1}\left(S_{t}^\top +B_{t}^\top P_{t}+D_{t}^\top P_{t}C_{t}+\int_{\Theta} F_{t,\theta}^\top P_{t}E_{t,\theta}\nu(d\theta)\right)(X_{t}^*-\mathbb E[X_{t}^*])\\&-\Sigma_{1t}^{-1}\left[(S_{t}+\bar S_{t})^\top+(B_{t}+\bar B_{t})^\top \Pi_{t}+(D_{t}+\bar D_{t})^\top P_{t}(C_{t}+\bar C_{t})\right.\\&\left.+\int_{\Theta} (F_{t,\theta}+\bar F_{t,\theta})^\top P_{t}(E_{t,\theta}+\bar E_{t,\theta})\nu(d\theta)\right]\mathbb E[X_{t}^*],\\
dX_{t}^*=&\ \Big\{A_{t}X_{t}^*+\bar A_{t}\mathbb E[X_{t}^*]+B_{t}u_{t}^*+\bar B_{t}\mathbb E[u_{t}^*]\Big\}dt
\\&+\Big\{C_{t}X_{t}^*+\bar C_{t}\mathbb E[X_{t}^*]+D_{t}u_{t}^*+\bar D_{t}\mathbb E[u_{t}^*]\Big\}dW_{t}
\\&+\int_{\Theta}\Big\{E_{t,\theta} X_{t-}^*+\bar E_{t,\theta}\mathbb E[ X_{t-}^*]+F_{t,\theta}u_{t}^*+\bar F_{t,\theta}\mathbb E[ u_{t}^*]\Big \}\widetilde N(dt,d\theta),\\
X_{0}^*=&\ x.
\end{aligned}
\right.
\end{equation*}
Defining
\begin{equation}\label{decouple MF}
\left\{\begin{aligned}
Y_{t}^*=&\ P_{t}(X_{t}^*-\mathbb E[X_{t}^*])+\Pi_{t}\mathbb E[X_{t}^*],
\\Z_{t}^*=&\ \left[P_{t}C_{t}-P_{t}D_{t}\Sigma_{0t}^{-1}\left(S_{t}^\top +B_{t}^\top P_{t}+D_{t}^\top P_{t}C_{t}+\int_{\Theta} F_{t,\theta}^\top P_{t} E_{t,\theta}\nu(d\theta)\right)\right](X_{t}^*-\mathbb E[X_{t}^*])
\\&+\left\{P_{t}(C_{t}+\bar C_{t})-P_{t}(D_{t}+\bar D_{t})\Sigma_{1t}^{-1}\left[(S_{t}+\bar S_{t})^\top+(B_{t}+\bar B_{t})^\top \Pi_{t}\right.\right.
\\&\left.\left.+(D_{t}+\bar D_{t})^\top P_{t}(C_{t}+\bar C_{t})+\int_{\Theta} (F_{t,\theta}+\bar F_{t,\theta})^\top P_{t}(E_{t,\theta}+\bar E_{t,\theta})\nu(d\theta)\right]\right\}\mathbb E[X_t^*],
\\
r_{t,\theta}^*=&\ \left[P_{t}E_{t,\theta}-P_{t}F_{t,\theta}\Sigma_{0t}^{-1}\left(S_{t}^\top +B_{t}^\top P_{t}+D_{t}^\top P_{t}C_{t}+\int_{\Theta} F_{t,\theta}^\top P_{t} E_{t,\theta}\nu(d\theta)\right)\right](X_{t}^*-\mathbb E[X_{t}^*])
\\&+\left\{P_{t}(E_{t,\theta}+\bar E_{t,\theta})-P_{t}(F_{t,\theta}+\bar F_{t,\theta})\Sigma_{1t}^{-1}\left[(S_{t}+\bar S_{t})^\top+(B_{t}+\bar B_{t})^\top \Pi_{t}\right.\right.
\\&+\left.\left.+(D_{t}+\bar D_{t})^\top P_{t}(C_{t}+\bar C_{t})\int_{\Theta} (F_{t,\theta}+\bar F_{t,\theta})^\top P_{t}(E_{t,\theta}+\bar E_{t,\theta})\nu(d\theta)\right]\right\}\mathbb E[\tilde X_{t}^*],
\end{aligned}
\right.
\end{equation}
the $5$-tuple $(X^*, u^*, Y^*,Z^*,r^*)$ is the unique  solution to  MF-FBSDEJ \eqref{H   system}. Moreover,
\begin{equation*}
\inf_{u\in \mathcal U[0,T]}J[u]=\frac{1}{2}\langle\Pi_{0}x,x\rangle,\ \ \ \forall x\in \mathbb R^n.
\end{equation*}
\end{thm}
\begin{proof}
We only need to prove  the $5$-tuple $(X^*,u^*, Y^*,Z^*,r^*)$ is  the unique solution to   MF-FBSDEJ \eqref{H   system}. By the linear transformation introduced in Lemma \ref{control system lemma} and the
representation of optimal control $\tilde u^*$ for Problem NC in Lemma \ref{Yong}, we get
\begin{equation*}
\begin{aligned}
u_{t}^*=&\ u_{t}^*-\mathbb E[u_{t}^*]+\mathbb E[u_{t}^*]
\\=&-R_{t}^{-1}S_{t}^\top(\tilde X_{t}^*-\mathbb E[\tilde X_{t}^*])+\tilde u_{t}^*-\mathbb E[\tilde u_{t}^*]-(R_{t}+\bar R_{t})^{-1}(S_{t}+\bar S_{t})^\top\mathbb E[\tilde X_{t}^*]+\mathbb E[\tilde u_{t}^*]
\\=&-R_{t}^{-1}S_{t}^\top( X_{t}^*-\mathbb E[ X_{t}^*])-\Sigma_{0t}^{-1}\left(B_{t}^\top P_{t}+D_{t}^\top P_{t}C_{1t}\right.\\&\left.+\int_{\Theta} F_{t,\theta}^\top P_{t}E_{1t,\theta}\nu(d\theta)\right)(X_{t}^*-\mathbb E[X_{t}^*])
-(R_{t}+\bar R_{t})^{-1}(S_{t}+\bar S_{t})^\top\mathbb E[X_{t}^*]
\\&-\Sigma_{1t}^{-1}\left[(B_{t}+\bar B_{t})^\top \Pi_{t}+(D_{t}+\bar D_{t})^\top P_{t}(C_{1t}+\bar C_{1t})\right.
\\&\left.+\int_{\Theta} (F_{t,\theta}+\bar F_{t,\theta})^\top P_{t}(E_{1t,\theta}+\bar E_{1t,\theta})\nu(d\theta)\right]\mathbb E[X_{t}^*]
\\=&-\Sigma_{0t}^{-1}\left[\Sigma_{0t}R_{t}^{-1}S_{t}^\top+B_{t}^\top P_{t}+D_{t}^\top P_{t}(C_{t}-D_{t}R_{t}^{-1}S_{t}^\top)\right.
\\&\left.+\int_{\Theta} F_{t,\theta}^\top P_{t} (E_{t,\theta}-F_{t,\theta}R_{t}^{-1}S_{t}^\top)\nu(d\theta)\right](X_{t}^*-\mathbb E[X_{t}^*])
\\&-\Sigma_{1t}^{-1} \left\{\Sigma_{1t}(R_{t}+\bar R_{t})^{-1}(S_{t}+\bar S_{t})^\top+(B_{t}+\bar B_{t})^\top \Pi_{t}\right.\\&+(D_{t}+\bar D_{t})^\top P_{t}\left[C_{t}+\bar C_{t}-(D_{t}+\bar D_{t})(R_{t}+\bar R_{t})^{-1}(S_{t}+\bar S_{t})^\top\right]
\\&\left.+\int_{\Theta} (F_{t,\theta}+\bar F_{t,\theta})^\top P_{t}\left[E_{t,\theta}+\bar E_{t,\theta}-(F_{t,\theta}+\bar F_{t,\theta})(R_{t}+\bar R_{t})^{-1}(S_{t}+\bar S_{t})^\top\right]\nu(d\theta)\right\}\mathbb E[X_{t}^*]
\\=&-\Sigma_{0t}^{-1}\left[S_{t}^\top +B_{t}^\top P_{t}+D_{t}^\top P_{t}C_{t}+\int_{\Theta} F_{t,\theta}^\top P_{t}E_{t,\theta}\nu(d\theta)\right](X_{t}^*-\mathbb E[X_{t}^*])
\\&-\Sigma_{1t}^{-1}\left[(S_{t}+\bar S_{t})^\top+(B_{t}+\bar B_{t})^\top \Pi_{t}+(D_{t}+\bar D_{t})^\top P_{t}(C_{t}+\bar C_{t})\right.\\&\left.+\int_{\Theta} (F_{t,\theta}+\bar F_{t,\theta})^\top P_{t}(E_{t,\theta}+\bar E_{t,\theta})\nu(d\theta)\right]\mathbb E[X_{t}^*].
\end{aligned}
\end{equation*}
Using Lemma \ref{H system lemma} and the representation of $\tilde Y^*, \tilde Z^*,\tilde r^*$ in Lemma \ref{Yong}, we have
\begin{equation*}
\begin{aligned}
Y_{t}^*=&\ \tilde Y_{t}^*\\=&P_{t}(\tilde X_{t}^*-\mathbb E[\tilde X_{t}^*])+\Pi_{t}\mathbb E[\tilde X_{t}^*]\\=&P_{t}( X_{t}^*-\mathbb E[ X_{t}^*])+\Pi_{t}\mathbb E[ X_{t}^*],
\\
Z_{t}^*=&\ \tilde Z_{t}^*
\\=&\left[P_{t}C_{1t}-P_{t}D_{t}\Sigma_{0t}^{-1}\left(B_{t}^\top P_{t}+D_{t}^\top P_{t}C_{1t}+\int_{\Theta} F_{t,\theta}^\top P_{t}E_{1t,\theta}\nu(d\theta)\right)\right](\tilde X_{t}^*-\mathbb E[\tilde X_{t}^*])
\\&+\left\{P_{t}(C_{1t}+\bar C_{1t})-P_{t}(D_{t}+\bar D_{t})\Sigma_{1t}^{-1}\left[(B_{t}+\bar B_{t})^\top \Pi_{t}+(D_{t}+\bar D_{t})^\top P_{t}(C_{1t}+\bar C_{1t})\right.\right.
\\&\left.\left.+\int_{\Theta} (F_{t,\theta}+\bar F_{t,\theta})^\top P_{t}(E_{1t,\theta}+\bar E_{1t,\theta})\nu(d\theta)\right]\right\}\mathbb E[\tilde X_{t}^*]
\\=&\left\{P_{t}(C_{t}-D_{t}R_{t}^{-1}S_{t}^\top)-P_{t}D_{t}\Sigma_{0t}^{-1}\left[B_{t}^\top P_{t}+D_{t}^\top P_{t}(C_{t}-D_{t}R_{t}^{-1}S_{t}^\top)\right.\right.
\\&\left.\left.+\int_{\Theta} F_{t,\theta}^\top P_{t}(E_{t,\theta}-F_{t,\theta}R_{t}^{-1}S_{t}^\top)\nu(d\theta)\right]\right\}( X_{t}^*-\mathbb E[X_{t}^*])
\\&+\left\{P_{t}\left[C_{t}+\bar C_{t}-(D_{t}+\bar D_{t})(R_{t}+\bar R_{t})^{-1}(S_{t}+\bar S_{t})^\top\right]-P_{t}(D_{t}+\bar D_{t})\Sigma_{1t}^{-1}\left[(B_{t}+\bar B_{t})^\top \Pi_{t}\right.\right.
\\&+(D_{t}+\bar D_{t})^\top P_{t}\left(C_{t}+\bar C_{t}-\left(D_{t}+\bar D_{t}\right)\left(R_{t}+\bar R_{t}\right)^{-1}\left(S_{t}+\bar S_{t}\right)^\top\right)
\\&\left.\left.+\int_{\Theta} \left(F_{t,\theta}+\bar F_{t,\theta}\right)^\top P_{t}\left(E_{t,\theta}+\bar E_{t,\theta}-\left(F_{t,\theta}+\bar F_{t,\theta}\right)\left(R_{t}+\bar R_{t}\right)^{-1}\left(S_{t}+\bar S_{t}\right)^\top\right)\nu(d\theta)\right]\right\}\mathbb E[X_{t}^*]
\\=&\left[P_{t}C_{t}-P_{t}D_{t}\Sigma_{0t}^{-1}\left(S_{t}^\top+B_{t}^\top P_{t}+D_{t}^\top P_{t}C_{t}+\int_{\Theta} F_{t,\theta}^\top P_{t} E_{t,\theta}\nu(d\theta)\right)\right](X_{t}^*-\mathbb E[X_{t}^*])
\\&+\left\{P_{t}(C_{t}+\bar C_{t})-P_{t}(D_{t}+\bar D_{t})\Sigma_{1t}^{-1}\left[(S_{t}+\bar S_{t})^\top+(B_{t}+\bar B_{t})^\top \Pi_{t}
\right.\right.\\&\left.\left.+(D_{t}+\bar D_{t})^\top P_{t}(C_{t}+\bar C_{t})+\int_{\Theta} (F_{t,\theta}+\bar F_{t,\theta})^\top P_{t}(E_{t,\theta}+\bar E_{t,\theta})\nu(d\theta)\right]\right\}\mathbb E[X_{t}^*].
\end{aligned}
\end{equation*}
We can also obtain that $r_{t,\theta}^*$ satisfies \eqref{decouple MF} similarly. Then the proof is completed.

\end{proof}

\section{Indefinite MF-LQJ problem}
In the case that Assumption (S) does not hold
true, it is possible that Problem MF is well-posed and an optimal pair exists.
In this section, we will apply the  equivalent cost functional method  to deal with
Problem MF without Assumption (S).
\begin{definition}
For a given controlled system, if there exist two cost functionals $J$ and $\bar J$ satisfying: for any admissible controls $u^1$ and $u^2$, $J[u^1]<J[u^2]$ if and only if $\bar J[u^1]<\bar J[u^2]$, then we say $J$ is equivalent to $\bar J$.
\end{definition}
\begin{rmk}
The following two statements are equivalent:
\begin{itemize}\item [1.]Cost functional $J$ is equivalent to $\bar J$;
 \item [2.] For any admissible controls $u^1$ and $u^2$,
    \begin{itemize}\item[(a).] $J[u^1]<J[u^2]$ if and only if $\bar J[u^1]<\bar J[u^2]$;
    \item[(b).]  $J[u^1]=J[u^2]$ if and only if $\bar J[u^1]=\bar J[u^2]$;
    \item[(c).]  $J[u^1]>J[u^2]$ if and only if $\bar J[u^1]>\bar J[u^2]$.
\end{itemize}\end{itemize}
\end{rmk}

We denote

$\Phi=\big\{\varphi\in L^\infty(0,T;\mathbb S^n)| \varphi$ is  a deterministic continuous differential function$\big\}$.

For any $H,K\in \Phi$, we define $J^{HK}[u]=J[u]-\frac{1}{2}\langle K_{0}x,x\rangle$.  Applying It\^o formula to $\langle H_{t}(X_{t}-\mathbb E[X_{t}]), X_{t}-\mathbb E[X_{t}]\rangle+ \langle K_{t}\mathbb E[X_{t}],\mathbb E[X_{t}]\rangle$ on the interval $[0,T]$, we derive
\begin{equation*}\label{cost HK function}
\begin{aligned}
&J^{HK}[u]\\=&
\ \frac{1}{2}\mathbb E\left\{\langle G^{HK}(X_{T}-\mathbb E[X_{T}]), X_{T}-\mathbb E[X_{T}]\rangle+\langle (G^{HK}+\bar G^{HK})\mathbb E[X_{T}],\mathbb E[X_{T}]\rangle\right.\\&+\int_0^T\left\langle\left(\begin{array}{cc}
Q_{t}^{HK}&{S_{t}^{HK}}^\top\\ S_{t}^{HK}&R_{t}^{HK}
\end{array}\right)\left(\begin{array}{c}
X_{t}-\mathbb{E}[X_{t}]\\ u_{t}-\mathbb{E}[u_{t}]
\end{array}\right),\left(\begin{array}{c}
X_{t}-\mathbb{E}[X_{t}]\\u_{t}-\mathbb{E}[u_{t}]
\end{array}\right)\right\rangle dt\\&\left.+\int_0^T\Bigg\langle\left(\begin{array}{cc}
Q_{t}^{HK}+\bar Q_{t}^{HK}&(S_{t}^{HK}+\bar S_{t}^{HK})^\top\\S_{t}^{HK}+\bar S_{t}^{HK}&R_{t}^{HK}+\bar R_{t}^{HK}
\end{array}\right)\left(\begin{array}{c}
\mathbb E[X_{t}]\\\mathbb E[u_{t}]
\end{array}\right),\left(\begin{array}{c}
\mathbb E[X_{t}]\\\mathbb E[u_{t}]
\end{array}\right)\Bigg\rangle dt\right\},
\end{aligned}
\end{equation*}
where
\begin{equation}\label{cost HK}
\begin{aligned}
&Q_{t}^{HK}= Q_{t}+\dot{H_{t}}+H_{t}A_{t}+A_{t}^\top H_{t}+C_{t}^\top H_{t}C_{t}+\int_{\Theta} E_{t,\theta}^\top H_{t}E_{t,\theta}\nu(d\theta),\\&
S_{t}^{HK}= S_{t}+H_{t}B_{t}+C_{t}^\top H_{t}D_{t}+\int_{\Theta} E_{t,\theta}^\top H_{t} F_{t,\theta}\nu(d\theta),\\&
R_{t}^{HK}= R_{t}+D_{t}^\top H_{t}D_{t}+\int_{\Theta} F_{t,\theta}^\top H_{t}F_{t,\theta}\nu(d\theta),\\&
G^{HK}= G-H_{T},\\&
Q_{t}^{HK}+\bar Q_{t}^{HK}= Q_{t}+\bar Q_{t}+\dot{K_{t}}+K_{t}(A_{t}+\bar A_{t})+(A_{t}+\bar A_{t})^\top K_{t}\\&\quad+(C_{t}+\bar C_{t})^\top H_{t}(C_{t}+\bar C_{t})+\int_{\Theta} (E_{t,\theta}+\bar E_{t,\theta})^\top H_{t}(E_{t,\theta}+\bar E_{t,\theta})\nu(d\theta),\\&
S_{t}^{HK}+\bar S_{t}^{HK}= S_{t}+\bar S_{t}+K_{t}(B_{t}+\bar B_{t})+(C_{t}+\bar C_{t})^\top H_{t}(D_{t}+\bar D_{t})\\&\quad+\int_{\Theta} (E_{t,\theta}+\bar E_{t,\theta})^\top H_{t}(F_{t,\theta}+\bar F_{t,\theta})\nu(d\theta),\\&
R_{t}^{HK}+\bar R_{t}^{HK}= R_{t}+\bar R_{t}+(D_{t}+\bar D_{t})^\top H_{t}(D_{t}+\bar D_{t})+\int_{\Theta} (F_{t,\theta}+\bar F_{t,\theta})^\top H_{t}(F_{t,\theta}+\bar F_{t,\theta})\nu(d\theta),\\&
G^{HK}+\bar G^{HK}= G+\bar G-K_{T}.
\end{aligned}
\end{equation}
Since $J^{HK}[u]$ and $J[u]$ differ by only a constant $-\frac{1}{2}\langle K_{0}x, x\rangle$, they are equivalent. In other words, we get a family of equivalent cost functionals $\{J^{HK}[u]\}$, which includes the original cost functional $J[u]$ when  $H\equiv 0, \ K\equiv 0$. For the sake
of convenience, given any $H, K\in \Phi$, we call  Problem MF ``Problem $J^{HK}$" if $J[u]$ is replaced by $J^{HK}[u]$.

We write down   stochastic  Hamiltonian systems corresponding to Problem $J^{HK}$
\begin{equation}\label{H   system1}
\left\{\begin{aligned}&dX_{t}^{HK}=\left\{A_{t}X_{t}^{HK}+\bar A_{t}\mathbb E[X_{t}^{HK}]+B_{t}u_{t}^{HK}+\bar B_{t}\mathbb E[u_{t}^{HK}]\right\}dt
\\&\quad+\left\{C_{t}X_{t}^{HK}+\bar C_{t}\mathbb E[X_{t}^{HK}]+D_{t}u_{t}^{HK}+\bar D_{t}\mathbb E[u_{t}^{HK}]\right\}dW_{t}
\\&\quad+\int_{\Theta}\left\{E_{t,\theta}X_{t-}^{HK}+\bar E_{t,\theta}\mathbb E[X_{t-}^{HK}]+F_{t,\theta}u_{t}^{HK}+\bar F_{t,\theta}\mathbb E[u_{t}^{HK}]\right\}\widetilde N(dt,d\theta),
\\&dY_{t}^{HK}=-\left\{A_{t}^\top Y_{t}^{HK}+\bar A_{t}^\top\mathbb E[Y_{t}^{HK}]+C_{t}^\top Z_{t}^{HK}+\bar C_{t}^\top \mathbb E[Z_{t}^{HK}]\right.\\&\quad+\int_{\Theta} E_{t,\theta}^\top r_{t,\theta}^{HK}\nu(d\theta)+\int_{\Theta} \bar E_{t,\theta}^\top \mathbb E[r_{t,\theta}^{HK}]\nu(d\theta)+Q_{t}^{HK}X_{t}^{HK}+\bar Q_{t}^{HK}\mathbb E[X_{t}^{HK}]\\&\quad\left.+S_{t}^{HK}u_{t}^{HK}+\bar S_{t}^{HK}\mathbb E[u_{t}^{HK}]\right\}dt+Z_{t}^{HK}dW+\int_{\Theta} r_{t,\theta}^{HK}\widetilde{N}(dt,d\theta), \\&X_{0}^{HK}=x, \qquad Y_{T}^{HK}=G^{HK}X_{T}^{HK}+\bar G^{HK}\mathbb E[X_{T}^{HK}],\\&{S_{t}^{HK}}^\top X_{t-}^{HK}+\bar{S}_{t}^{{HK}^\top}\mathbb E[X_{t-}^{HK}]+B_{t}^\top Y_{t-}^{HK}+\bar B_{t}^\top\mathbb E[Y_{t-}^{HK}]+D_{t}^\top Z_{t}^{HK}+\bar D_{t}^\top\mathbb E[Z_{t}^{HK}]\\&+\int_{\Theta} F_{t,\theta}^\top r_{t,\theta}^{HK}\nu(d\theta)+\int_{\Theta} \bar F_{t,\theta}^\top \mathbb E[r_{t,\theta}^{HK}]\nu(d\theta)+R_{t}^{HK}u_{t}^{HK}+\bar R_{t}^{HK}\mathbb E[u_{t}^{HK}]=0.
\end{aligned}
\right.
\end{equation}
We note that \eqref{H   system} coincides with \eqref{H   system1} while $H\equiv 0, K\equiv 0$. The following lemma shows that there  exists an equivalent relationship among Hamiltonian systems  \eqref{H   system1}.

\begin{lem}\label{H system Indefinite}
For all $H,K\in \Phi$, the existence and uniqueness of solutions  to Hamiltonian systems  \eqref{H   system1} are equivalent.
\end{lem}
\begin{proof}
We only need to prove that for all $H,K\in \Phi$, the existence and uniqueness of solutions to Hamiltonian systems \eqref{H   system1} are equivalent to \eqref{H   system}. For any given $H,K\in \Phi$, if $(X^{HK}, u^{HK}, Y^{HK}, Z^{HK}, r^{HK})$ is a solution of \eqref{H   system1}, we define
\begin{equation*}
\begin{aligned}
&X_{t}= X_{t}^{HK},\\&
u_{t}= u_{t}^{HK},\\&
Y_{t}= Y_{t}^{HK}+H_{t}(X_{t}^{HK}-\mathbb E[X_{t}^{HK}])+K_{t}\mathbb E[X_{t}^{HK}],\\&
Z_{t}= Z_{t}^{HK}+H_{t}(C_{t}X_{t}^{HK}+\bar C_{t}\mathbb E[X_{t}^{HK}]+D_{t}u_{t}^{HK}+\bar D_{t}\mathbb E[u_{t}^{HK}]),\\&
r_{t,\theta}= r_{t,\theta}^{HK}+H_{t}(E_{t,\theta}X_{t-}^{HK}+\bar E_{t,\theta}\mathbb E[X_{t-}^{HK}]+F_{t,\theta}u_{t}^{HK}+\bar F_{t,\theta}\mathbb E[u_{t}^{HK}]).
\end{aligned}
\end{equation*}
Then   $(X, u, Y, Z, r)$ is a solution of \eqref{H   system}. Thus we complete the proof.

\end{proof}

We write down    Riccati equations related to Problem $J^{HK}$:

\begin{equation}\label{P riccati1}
\left\{\begin{aligned}&\dot{P_{t}}^{HK}+P_{t}^{HK}A_{t}+A_{t}^\top P_{t}^{HK}+C_{t}^\top P_{t}^{HK}C_{t}+\int_{\Theta} E_{t,\theta}^\top P_{t}^{HK} E_{t,\theta}\nu(d\theta)+Q_{t}^{HK}\\&-\left( S_{t}^{HK}+P_{t}^{HK}B_{t}+C_{t}^\top P_{t}^{HK}D_{t}+\int_{\Theta} E_{t,\theta}^\top P_{t}^{HK} F_{t,\theta}\nu(d\theta)\right){\Sigma_{0t}^{HK}}^{-1}\\&\cdot\left({S_{t}^{HK}}^\top+B_{t}^\top P_{t}^{HK}+D_{t}^\top P_{t}^{HK}C_{t}+\int_{\Theta} F_{t,\theta}^\top P_{t}^{HK} E_{t,\theta}\nu(d\theta)\right)=0,\\ &P^{HK}_{T}=G^{HK},
\end{aligned}
\right.
\end{equation}

and

\begin{equation}\label{Pi riccati1}
\left\{\begin{aligned}
&\dot{\Pi_{t}}^{HK}+\Pi_{t}^{HK}(A_{t}+\bar A_{t})+(A_{t}+\bar A_{t})^\top \Pi_{t}^{HK}+(C_{t}+\bar C_{t})^\top P_{t}^{HK}(C_{t}+\bar C_{t})
\\&+\int_{\Theta} (E_{t,\theta}+\bar E_{t,\theta})^\top P_{t}^{HK} (E_{t,\theta}+\bar E_{t,\theta})\nu(d\theta)+Q_{t}^{HK}+\bar Q_{t}^{HK}
\\&-\left[(S_{t}^{HK}+\bar S_{t}^{HK})+\Pi_{t}^{HK}(B_{t}+\bar B_{t})+(C_{t}+\bar C_{t})^\top P_{t}^{HK}(D_{t}+\bar D_{t})\right.
\\&+\left.\int_{\Theta} (E_{t,\theta}+\bar E_{t,\theta})^\top P_{t}^{HK} (F_{t,\theta}+\bar  F_{t,\theta})\nu(d\theta)\right]{\Sigma_{1t}^{HK}}^{-1}
\\&\cdot\left[(S_{t}^{HK}+\bar S_{t}^{HK})^\top+(B_{t}+\bar B_{t})^\top \Pi_{t}^{HK}+(D_{t}+\bar D_{t})^\top P_{t}^{HK}(C_{t}+\bar C_{t})\right.
\\&+\left.\int_{\Theta} (F_{t,\theta}+\bar F_{t,\theta})^\top P_{t}^{HK} (E_{t,\theta}+\bar E_{t,\theta})\nu(d\theta)\right]=0,\\ &\Pi^{HK}_{T}=G^{HK}+\bar G^{HK},
\end{aligned}
\right.
\end{equation}

where
\begin{equation*}
\begin{aligned}
&\Sigma_{0t}^{HK}=R_{t}^{HK}+D_{t}^\top P_{t}^{HK}D_{t}+\int_{\Theta} F_{t,\theta}^\top P_{t}^{HK} F_{t,\theta}\nu(d\theta),
\\&\Sigma_{1t}^{HK}=R_{t}^{HK}+\bar R_{t}^{HK}+(D_{t}+\bar D_{t})^\top P_{t}^{HK}(D_{t}+\bar D_{t})+\int_{\Theta} (F_{t,\theta}+\bar F_{t,\theta})^\top P_{t}^{HK} (F_{t,\theta}+\bar F_{t,\theta})\nu(d\theta).
\end{aligned}
\end{equation*}
\begin{lem}\label{Riccati  Indefiite}
For any $H,K\in \Phi$, the existence and uniqueness of solutions to Riccati equations associated with {\it Problem} $J^{HK}$ are equivalent.
\end{lem}
\begin{proof}
We only need to prove that for all $H,K\in \Phi$, the existence and uniqueness
of solutions to Riccati equations \eqref{P riccati1} and \eqref{Pi riccati1} are equivalent to Riccati equations \eqref{P riccati} and \eqref{Pi riccati}.
For any given $H,K\in \Phi$, if $P^{HK}, \Pi^{HK}$ are solutions of Riccati equations \eqref{P riccati1} and \eqref{Pi riccati1},
respectively, then we define
\begin{equation*}
\begin{aligned}
P_{t}=P_{t}^{HK}+H_{t},\\
\Pi_{t}=\Pi_{t}^{HK}+K_{t}.
\end{aligned}
\end{equation*}
Through a straightforward calculation, we know that $P,\Pi$ are solutions of Riccati equations
\eqref{P riccati} and \eqref{Pi riccati}, respectively. Thus we complete the proof.
 \end{proof}

We now give two theorems which are our main results in this section.
\begin{thm}\label{main results Th}
If there exist  $\bar H,\bar K\in \Phi$ such that $(Q_t^{\bar H\bar K},\ S_t^{\bar H\bar K},\ R_t^{\bar H\bar K},\ G^{\bar H\bar K})$,   $(\bar Q_t^{\bar H\bar K},\ \bar S_t^{\bar H\bar K},\ \bar R_t^{\bar H\bar K},\ \bar G^{\bar H\bar K})$ satisfy Assumption (S), then    Riccati equations \eqref{P riccati} and \eqref{Pi riccati} admit  unique solutions $P,\Pi$,  respectively. Further, the
optimal pair $(X^*, u^*)$ of {\it Problem MF} satisfies
\begin{equation}\label{equipair}
\left\{\begin{aligned}
u_{t}^*=&-\Sigma_{0t}^{-1}\left(S_{t}^\top+B_{t}^\top P_{t}+D_{t}^\top P_{t}C_{t}+\int_{\Theta} F_{t,\theta}^\top P_{t}^{HK} E_{t,\theta}\nu(d\theta)\right)(X_{t}^*-\mathbb E[X_{t}^*])
\\&-\Sigma_{1t}^{-1}\left[(S_{t}+\bar S_{t})^\top+(B_{t}+\bar B_{t})^\top \Pi_{t}+(D_{t}+\bar D_{t})^\top P_{t}(C_{t}+\bar C_{t})\right.
\\&\left.+\int_{\Theta} (F_{t,\theta}+\bar F_{t,\theta}) ^\top P_{t}(E_{t,\theta}+\bar E_{t,\theta})\nu(d\theta)\right]\mathbb E[X_{t}^*],
\\
dX_{t}^*=&\Big(A_{t}X_{t}^*+\bar A_{t}\mathbb E[X_{t}^*]+B_{t}u_{t}^*+\bar B_{t}\mathbb E[u_{t}^*]\Big)dt
\\&+\Big(C_{t}X_{t}^*+\bar C_{t}\mathbb E[X_{t}^*]+D_{t}u_{t}^*+\bar D_{t}\mathbb E[u_{t}^*]\Big)dW_{t}
\\&+\int_{\Theta}\Big\{E_{t,\theta}X_{t-}^*+\bar E_{t,\theta}\mathbb E[X_{t-}^*]+F_{t,\theta}u_{t}^*+\bar F_{t,\theta}\mathbb E[u_{t}^*]\Big \}\widetilde N(dt,d\theta),\\
X_{0}^*=&\ x.
\end{aligned}
\right.
\end{equation}
Defining
\begin{equation}\label{equiadjoint}
\left\{\begin{aligned}
Y_{t}^*=&\ P_{t}(X_{t}^*-\mathbb E[X_{t}^*])+\Pi_{t}\mathbb E[X_{t}^*],
\\Z_{t}^*=&\ \left[P_{t}C_{t}-P_{t}D_{t}\Sigma_{0t}^{-1}\left(S_{t}^\top +B_{t}^\top P_{t}+D_{t}^\top P_{t}C_{t}+\int_{\Theta} F_{t,\theta}^\top P_{t}  E_{t,\theta}\nu(d\theta)\right)\right](X_{t}^*-\mathbb E[X_{t}^*])
\\&+\left\{P_{t}(C_{t}+\bar C_{t})-P_{t}(D_{t}+\bar D_{t})\Sigma_{1t}^{-1}\left[(S_{t}+\bar S_{t})^\top+(B_{t}+\bar B_{t})^\top \Pi_{t}
\right.\right.
\\&\left.\left.+(D_{t}+\bar D_{t})^\top P_{t}(C_{t}+\bar C_{t})+\int_{\Theta} (F_{t,\theta}+\bar F_{t,\theta})^\top P_{t}(E_{t,\theta}+\bar E_{t,\theta})\nu(d\theta)\right]\right\}\mathbb E[X_{t}^*],
\\r_{t,\theta}^*=&\ \left[P_{t}E_{t,\theta}-P_{t}F_{t,\theta}\Sigma_{0t}^{-1}\left(S_{t}^\top+B_{t}^\top P_{t}+D_{t}^\top P_{t}C_{t}+\int_{\Theta} F_{t,\theta}^\top P_{t}E_{t,\theta}\nu(d\theta)\right)\right](X_{t}^*-\mathbb E[X_{t}^*])
\\&+\Big\{P_{t}(E_{t,\theta}+\bar E_{t,\theta})-P_{t}(F_{t,\theta}+\bar F_{t,\theta})\Sigma_{1t}^{-1}\Big[(S_{t}+\bar S_{t})^\top+(B_{t}+\bar B_{t})^\top \Pi_{t}
\\&+(D_{t}+\bar D_{t})^\top P_{t}(C_{t}+\bar C_{t})+\int_{\Theta} (F_{t,\theta}+\bar F_{t,\theta})^\top P_{t}(E_{t,\theta}+\bar E_{t,\theta})\nu(d\theta)\Big]\Big\}\mathbb E[\tilde X_{t}^*],
\end{aligned}
\right.
\end{equation}
the $5$-tuple $(X^*,u^*, Y^*,Z^*, r^*)$ is the unique solution to  MF-FBSDEJ \eqref{H   system}. Moreover,
\begin{equation*}
\inf_{u\in \mathcal U[0,T]}J[u]=\frac{1}{2}\langle\Pi_{0}x,x\rangle,\ \ \ \forall x\in \mathbb R^n.
\end{equation*}
\end{thm}
\begin{proof}
We now consider Problem $J^{\bar H\bar K}$. For $\bar H,\bar K\in \Phi$, according to Theorem \ref{MF STN},   Riccati equations \eqref{P riccati1} and \eqref{Pi riccati1} associated with Problem $J^{\bar H\bar K}$ admit unique solutions $P^{\bar H\bar K}, \Pi^{\bar H\bar K}$, respectively. By Lemma \ref{Riccati  Indefiite}, Riccati equations \eqref{P riccati} and \eqref{Pi riccati} admit  unique solutions $P,\Pi$, respectively.
Using Theorem \ref{MF STN}, further we  know that the
optimal pair $(X^{\bar H\bar K}, u^{\bar H\bar K})$ of Problem $J^{\bar H\bar K}$ satisfies
\begin{equation*}
\left\{\begin{aligned}
u_{t}^{\bar H\bar K}=
&-{\Sigma_{0t}^{\bar H\bar K}}^{-1}\left({S_{t}^{\bar H\bar K}}^\top +B_{t}^\top P_{t}^{\bar H\bar K}+D_{t}^\top P_{t}^{\bar H\bar K}C_{t}\right.\\&\left.+\int_{\Theta} F_{t,\theta}^\top P_{t}^{\bar H\bar K} E_{t,\theta}\nu(d\theta)\right)\left(X_{t}^{\bar H\bar K}-\mathbb E[ X_{t}^{\bar H\bar K}]\right)
\\&-{\Sigma_{1t}^{\bar H\bar K}}^{-1}\left[(S_{t}^{\bar H\bar K}+\bar S_{t}^{\bar H\bar K})^\top+(B_{t}+\bar B_{t})^\top \Pi_{t}^{\bar H\bar K}+(D_{t}+\bar D_{t})^\top P_{t}^{\bar H\bar K}(C_{t}+\bar C_{t})\right.\\&\left.+\int_{\Theta} (F_{t,\theta}+\bar F_{t,\theta})^\top P_{t}^{\bar H\bar K} (E_{t,\theta}+\bar E_{t,\theta})\nu(d\theta)\right]\mathbb E[X_{t}^{\bar H\bar K}],\\
dX_{t}^{\bar H\bar K}=&\ \left(A_{t}X_{t}^{\bar H\bar K}+\bar A_{t}\mathbb E[X_{t}^{\bar H\bar K}]+B_{t}u_{t}^{\bar H\bar K}+\bar B_{t}\mathbb E[u_{t}^{\bar H\bar K}]\right)dt\\&+\left(C_{t}X_{t}^{\bar H\bar K}+\bar C_{t}\mathbb E[ X_{t}^{\bar H\bar K}]+D_{t}u_{t}^{\bar H\bar K}+\bar D_{t}\mathbb E[ u_{t}^{\bar H\bar K}]\right)dW_{t}\\&+\int_{\Theta}\left\{E_{t,\theta}X_{t-}^{\bar H\bar K}+\bar E_{t,\theta}\mathbb E[X_{t-}^{\bar H\bar K}]+F_{t,\theta}u_{t}^{\bar H\bar K}+\bar F_{t,\theta}\mathbb E[u_{t}^{\bar H\bar K}]\right\}\widetilde N(dt,d\theta),\\
 X_{0}^{\bar H\bar K}=&\ x.
\end{aligned}
\right.
\end{equation*}
Defining
\begin{equation*}
\left\{\begin{aligned}
Y_{t}^{\bar H\bar K}=&\ P_{t}^{\bar H\bar K}\left(X_{t}^{\bar H\bar K}-\mathbb E[X_{t}^{\bar H\bar K}]\right)+\Pi_{t}^{\bar H\bar K}\mathbb E[X_{t}^{\bar H\bar K}],
\\ Z_{t,\theta}^{\bar H\bar K}=&\ \left[P_{t}^{\bar H\bar K}C_{t}-P_{t}^{\bar H\bar K}D_{t}{\Sigma_{0t}^{\bar H\bar K}}^{-1}\Big({S_{t}^{\bar H\bar K}}^\top+B_{t}^\top P_{t}^{\bar H\bar K}+D_{t}^\top P_{t}^{\bar H\bar K}C_{t}\right.
\\&\left.+\int_{\Theta}F_{t,\theta}^\top P_{t}^{\bar H\bar K}E_{t,\theta}\nu(d\theta)\Big)\right]\left(X_{t}^{\bar H\bar K}-\mathbb E[X_{t}^{\bar H\bar K}]\right)
\\&+\Big[P_{t}^{\bar H\bar K}(C_{t}+\bar C_{t})-P_{t}^{\bar H\bar K}(D_{t}+\bar D_{t}){\Sigma_{1t}^{\bar H\bar K}}^{-1}\Big((S_{t}^{\bar H\bar K}+\bar S_{t}^{\bar H\bar K})^\top+(B_{t}+\bar B_{t})^\top \Pi_{t}^{\bar H\bar K}
\\&+(D_{t}+\bar D_{t})^\top P_{t}^{\bar H\bar K}(C_{t}+\bar C_{t})+\int_{\Theta} (F_{t,\theta}+\bar F_{t,\theta})^\top P_{t}^{\bar H\bar K}(E_{t,\theta}+\bar E_{t,\theta})\nu(d\theta)\Big)\Big]\mathbb E[X^{\bar H\bar K}],
\\ r_{t,\theta}^{\bar H\bar K}=&\ \left[P_{t}^{\bar H\bar K}E_{t,\theta}-P_{t}^{\bar H\bar K}F_{t,\theta}{\Sigma_{0t}^{\bar H\bar K}}^{-1}\Big({S_{t}^{\bar H\bar K}}^\top +B_{t}^\top P_{t}^{\bar H\bar K}+D_{t}^\top P_{t}^{\bar H\bar K}C_{t}\right.\\&+\left.\int_{\Theta}F_{t,\theta}^\top P_{t}^{\bar H\bar K}E_{t,\theta}\nu(d\theta)\Big)\right]\left(X_{t}^{\bar H\bar K}-\mathbb E[X_{t}^{\bar H\bar K}]\right)\\&+\Big[P_{t}^{\bar H\bar K}(E_{t,\theta}+\bar E_{t,\theta})-P_{t}^{\bar H\bar K}(F_{t,\theta}+\bar F_{t,\theta}){\Sigma_{1t}^{\bar H\bar K}}^{-1}\Big((S_{t}^{\bar H\bar K}+\bar S_{t}^{\bar H\bar K})^\top+(B_{t}+\bar B_{t})^\top \Pi_{t}^{\bar H\bar K}\\&+(D_{t}+\bar D_{t})^\top P_{t}^{\bar H\bar K}(C_{t}+\bar C_{t})+\int_{\Theta} (F_{t,\theta}+\bar F_{t,\theta})^\top P_{t}^{\bar H\bar K}(E_{t,\theta}+\bar E_{t,\theta})\nu(d\theta)\Big)\Big]\mathbb E[X_{t}^{\bar H\bar K}],
\end{aligned}
\right.
\end{equation*}
the $5$-tuple $(X^{\bar H\bar K},u^{\bar H \bar K}, Y^{\bar H \bar K},Z^{\bar H \bar K}, r^{\bar H \bar K})$ is a solution to MF-FBSDEJ \eqref{H   system1} corresponding to Problem $J^{\bar H\bar K}$. Since $J^{\bar H\bar K}[u]$ is equivalent to  $J[u]$, then $(X^{\bar H\bar K},u^{\bar H\bar K})$ is also an  optimal pair of Problem MF. By Lemma \ref{H system Indefinite}, we define
\begin{equation*}
\begin{aligned}
&X_{t}^*= X_{t}^{\bar H\bar K},\\&
u_{t}^*= u_{t}^{\bar H\bar K},\\&
Y_{t}^*= Y_{t}^{\bar H\bar K}+\bar H_{t}(X_{t}^{\bar H\bar K}-\mathbb E[X_{t}^{\bar H\bar K}])+\bar K\mathbb E[X_{t}^{\bar H\bar K}],\\&
Z_{t}^*= Z_{t}^{\bar H\bar K}+\bar H_{t}(C_{t}X_{t}^{\bar H\bar K}+\bar C_{t}\mathbb E[X_{t}^{\bar H\bar K}]+D_{t}u_{t}^{\bar H\bar K}+\bar D_{t} \mathbb E[u_{t}^{\bar H\bar K}]),
\\&r_{t,\theta}^*= r_{t,\theta}^{\bar H\bar K}+\bar H_{t}(E_{t,\theta}X_{t-}^{\bar H\bar K}+\bar E_{t,\theta}\mathbb E[X_{t-}^{\bar H\bar K}]+F_{t,\theta}u_{t}^{\bar H\bar K}+\bar F_{t,\theta} \mathbb E[u_{t}^{\bar H\bar K}]),
\end{aligned}
\end{equation*}
then the $5$-tuple $(X^*,u^*, Y^*,Z^*, r^*)$ is a solution to MF-FBSDEJ \eqref{H   system}. Through a simple calculation, we know  that the optimal pair $(X^*,u^*)$ satisfies equation \eqref{equipair} and $(Y^*,Z^*, r^*)$ satisfies \eqref{equiadjoint}.
Moreover, it follows from  Theorem \ref{MF STN} that
\begin{equation*}
\inf_{u\in \mathcal U[0,T]}J^{\bar H\bar K}[u]=\frac{1}{2}\langle\Pi_{0}^{\bar H\bar K}x,x\rangle,\ \ \ \forall x\in \mathbb R^n.
\end{equation*}
Thus we have
\begin{equation*}
\begin{aligned}
\inf_{u\in \mathcal U[0,T]}J[u]&=\frac{1}{2}\langle\Pi_{0}^{\bar H\bar K}x,x\rangle+\frac{1}{2}\langle K_{0}x,x\rangle\\&=\frac{1}{2}\langle\Pi_{0}x,x\rangle,\ \ \ \forall x\in \mathbb R^n.
\end{aligned}
\end{equation*}
The proof is completed.

\end{proof}
\begin{thm}\label{Theorem R}
If   Riccati equations \eqref{P riccati} and \eqref{Pi riccati} admit  unique solutions $P, \Pi$, respectively,  and   $R_{t}+D^\top P_{t}D_{t}+\int_{\Theta} F_{t,\theta}^\top P_{t}F_{t,\theta}\nu(d\theta)>\alpha_1I$, $R_{t}+\bar R_{t}+(D_{t}+\bar D_{t})^\top P_{t}(D_{t}+\bar D_{t})+\int_{\Theta} (F_{t,\theta}+\bar F_{t,\theta})^\top P_{t}(F_{t,\theta}+\bar F_{t,\theta})\nu(d\theta)>\alpha_1I$ for some $\alpha_1>0$, then there exist  $\bar H,\bar K\in \Phi$ such that $(Q_t^{\bar H\bar K}, S_t^{\bar H\bar K}, R_t^{\bar H\bar K}, G^{\bar H\bar K})$ and $(\bar Q_t^{\bar H\bar K}, \bar S_t^{\bar H\bar K}, \bar R_t^{\bar H\bar K}, \bar G^{\bar H\bar K})$ satisfy Assumption (S).
\end{thm}
\begin{proof}We consider the equivalent cost functional $J^{P\Pi}[u]$. It is easy to verify
\begin{equation*}
\begin{aligned}
&Q_{t}^{P\Pi}= Q_{t}+\dot{P_{t}}+P_{t}A_{t}+A_{t}^\top P_{t}+C_{t}^\top P_{t}C_{t}+\int_{\Theta} E_{t,\theta}^\top P_{t}E_{t,\theta}\nu(d\theta),\\&
S_{t}^{P\Pi}= S_{t}+B_{t}^\top P_{t}+D_{t}^\top P_{t}C_{t}+\int_{\Theta} F_{t,\theta}^\top P_{t}E_{t,\theta}\nu(d\theta),\\&
R_{t}^{P\Pi}= R_{t}+D_{t}^\top P_{t}D_{t}+\int_{\Theta} F_{t,\theta}^\top P_{t}F_{t,\theta}\nu(d\theta),\\&
G^{P\Pi}= G-P_{T},\\&
Q_{t}^{P\Pi}+\bar Q_{t}^{P\Pi}= Q_{t}+\bar Q_{t}+\dot{\Pi_{t}}+\Pi_{t}(A_{t}+\bar A_{t})+(A_{t}+\bar A_{t})^\top \Pi_{t}\\&\quad+(C_{t}+\bar C_{t})^\top P_{t}(C_{t}+\bar C_{t})+\int_{\Theta} (E_{t,\theta}+\bar E_{t,\theta})^\top P_{t}(E_{t,\theta}+\bar E_{t,\theta})\nu(d\theta),\\&
S_{t}^{P\Pi}+\bar S_{t}^{P\Pi}= S_{t}+\bar S_{t}+\Pi_{t}(B_{t}+\bar B_{t})+(C_{t}+\bar C_{t})^\top P_{t}(D_{t}+\bar D_{t})\\&\quad+\int_{\Theta} (E_{t,\theta}+\bar E_{t,\theta})^\top P_{t}(F_{t,\theta}+\bar F_{t,\theta})\nu(d\theta),\\&
R_{t}^{P\Pi}+\bar R_{t}^{P\Pi}= R_{t}+\bar R_{t}+(D_{t}+\bar D_{t})^\top P_{t}(D_{t}+\bar D_{t})+\int_{\Theta} (F_{t,\theta}+\bar F_{t,\theta})^\top P_{t}(F_{t,\theta}+\bar F_{t,\theta})\nu(d\theta),\\&
G^{P\Pi}+\bar G^{P\Pi}= G+\bar G-\Pi_{T},
\end{aligned}
\end{equation*}
and $(Q_t^{P\Pi}, S_t^{P\Pi}, R_t^{P\Pi}, G^{P\Pi})$, $(\bar Q_t^{P\Pi}, \bar S_t^{P\Pi}, \bar R_t^{P\Pi}, \bar G^{P\Pi})$ satisfy Assumption (S).

\end{proof}

\section{Examples}
In this section, we present four illustrative examples, where Assumption (S) does not hold true for original optimal control problems. Example 5.1  shows that an optimal control  exists even though Assumption (S) does not hold true. In Example 5.2, it is difficult to prove the existence and uniqueness of  solutions  to  related Riccati equations. We use the  equivalent cost functional method to construct an MF-LQJ problem which
satisfies Assumption (S) first, and then we obtain an optimal control of the original stochastic control problem via  solutions of Riccati equations. We   also give the existence and uniqueness of solutions for a family of MF-FBSDEJs as a byproduct of our results.  With the in-depth study of Example 5.2, we apply our results to prove the existence and uniqueness of solution to an MF-FBSDEJ in Example 5.3, where existing methods in literature do not work.  In Example 5.4, we apply our results to solve an  asset-liability management problem and give some  numerical solutions.

{\it\textbf{Example 5.1}}: Consider a 1-dimensional controlled MF-SDEJ
\begin{equation}\label{example 2}
\left\{ \begin{aligned}
dX_{t}=&\left(X_{t}-\mathbb E[X_{t}]+u_{t}+\mathbb E[u_{t}]\right)dt+(2u_{t}-\mathbb E[u_{t}])dW_{t}+\int_{[1,+\infty)}e^{-\theta}\mathbb E[u_{t}]\widetilde N(dt,d\theta),\\
X_{0}=&\ x,
\end{aligned}\right.
\end{equation}
with a cost functional
\begin{equation}\label{cost functional example2}
\begin{aligned}
J[u]=&\ \frac{1}{2}\mathbb E\Bigg\{2|X(T)|^2-|\mathbb E[X(T)]|^2\\&+\int_0^1\Big(-3|X_{t}|^2+3|\mathbb E[X_{t}]|^2-4|u_{t}|^2+2|\mathbb E[u_{t}]|^2\Big)dt\Bigg\}.
\end{aligned}
\end{equation}
With the data, Assumption (S) does not hold. We write down the stochastic Hamiltonian system
\begin{equation}\label{H   system example 2}
\left\{\begin{aligned}&dX_{t}=\left(X_{t}-\mathbb E[X_{t}]+u_{t}+\mathbb E[u_{t}]\right)dt+(2u_{t}-\mathbb E[u_{t}])dW_{t}+\int_{[1,+\infty)}e^{-\theta}\mathbb E[u_{t}]\widetilde N(dt,d\theta),\\&dY_{t}=-\Big\{Y_{t}-\mathbb E[Y_{t}]-3X_{t}+3\mathbb E[X_{t}]\Big\}dt+Z_{t}dW+\int_{[1,+\infty)} r_{t,\theta}\widetilde N(dt,d\theta),\\&X_{0}=x,\ \ \ \  Y(T)=2X(T)-\mathbb E[X(T)],\\&-4u_{t}+2\mathbb E[u_{t}]+Y_{t-}+\mathbb E[Y_{t-}]+2Z_{t}-\mathbb E[Z_{t}]+\int_{[1,+\infty)} e^{-\theta}\mathbb E[r_{t,\theta}]\nu(d\theta)=0.
\end{aligned}
\right.
\end{equation}
The corresponding Riccati equations   are
\begin{equation*}
\left\{\begin{aligned}&\dot{P_{t}}+2P_{t}-3-\frac{P_{t}^2}{4P_{t}-4}=0,\\
&P_{T}=2,
\end{aligned}
\right.
\end{equation*}
\begin{equation*}
\left\{\begin{aligned}&\dot{\Pi_{t}}-\frac{4\Pi_{t}^2}{-2+P_{t}+\delta P_{t}}=0,\\
&\Pi_{T}=1,
\end{aligned}
\right.
\end{equation*}
where $\delta=\int_{[1,+\infty)}e^{-2\theta}\nu(d\theta)>0$.
Solving them, we get  $$P_{t}=2, \ \Pi_{t}=\frac{\delta}{2T-2t+\delta}.$$
Note that
\begin{equation*}
\left\{\begin{aligned}&-4+4P_{t}=4,\\
&-2+P_{t}+\delta P_{t}=2\delta.
\end{aligned}
\right.
\end{equation*}
It follows from Theorem \ref{Theorem R} that  the equivalent cost functional $J^{P\Pi}[u]$ satisfies Assumption (S). According to  Theorem \ref{main results Th}, the optimal pair $(X, u)$ satisfies
\begin{equation*}
\left\{\begin{aligned}u_{t}=&-\frac{1}{2}(X_{t}-\mathbb E[X_{t}])-\frac{1}{2T -2t+\delta}\mathbb E[X_{t}],\\dX_{t}=&\left(X_{t}-\mathbb E[X_{t}]+u_{t}+\mathbb E[u_{t}]\right)dt+(2u_{t}-\mathbb E[u_{t}])dW_{t}+\int_{[1, +\infty)} e^{-\theta} \mathbb E[u_{t}]\widetilde N(dt,d\theta),\\
X_{0}=&\ x.
\end{aligned}
\right.
\end{equation*}
Defining
\begin{equation*}
\left\{\begin{aligned}Y_{t}=&\ 2(X_{t}-\mathbb E[X_{t}])-\frac{\delta}{2T-2t+\delta}\mathbb E[X_{t}],\\Z_{t}=&\ -2(X_{t}-\mathbb E[X_{t}])+\frac{2}{2T-2t+\delta}\mathbb E[X_{t}],\\r_{t,\theta}=&\- \frac{2e^{-\theta}}{2T-2t+\delta}\mathbb E[X_{t}],
\end{aligned}
\right.
\end{equation*}
the $5$-tuple $(X,u, Y,Z,r)$ is  a solution to MF-FBSDEJ \eqref{H   system example 2}. Moreover,
\begin{equation*}
\inf_{u\in \mathcal U[0,T]}J[u]=\frac{\delta}{2(2T+\delta)}x^2,\ \ \ \forall x\in \mathbb R.
\end{equation*}

{\it\textbf{Example 5.2}}:
Consider a 1-dimensional controlled MF-SDEJ
\begin{equation}\label{example 1}
\left\{ \begin{aligned}
dX_{t}=&\ \big(2X_{t}-\mathbb E[X_{t}]+u_{t}\big)dt+2u_{t}dW_{t}+\int_{\Theta} (E_{t,\theta}X_{t-}+\bar E_{t,\theta}\mathbb E[X_{t-}])\widetilde N(dt,d\theta),\quad t\in{[0,T]},\\
X_{0}=&\ x,
\end{aligned}\right.
\end{equation}
with a cost functional
\begin{equation}\label{cost functional example1}
\begin{aligned}
J[u]=&\frac{1}{2}\mathbb E\Big\{\alpha X_{T}^2-(\alpha+1)\mathbb E[X_{T}]^2+\int_0^T\left(4\mathbb E[X_{t}]^2+4\mathbb E[X_{t}]\mathbb E[u_{t}]+R_{t}u_{t}^2+\bar R_{t}\mathbb E[u_{t}]^2\right)dt\Big\},
\end{aligned}
\end{equation}
where $$\alpha >\frac{1}{2} (T+1)^2, R_{t}=(t+1)^3-2(t+1)^2, \bar R_{t}=1-(t+1)^3.$$

Clearly, Assumption (S) does not hold. Now we introduce an equivalent cost functional $J^{H_0K_0}[u]$ satisfying Assumption (S). Recalling \eqref{cost HK}, we have
\begin{equation*}
\begin{aligned}
&Q_{t}^{HK}= \dot{H_{t}}+4H_{t}+\delta_{1t} H_{t},\
S_{t}^{HK}= H_{t},\\&
R_{t}^{HK}= (t+1)^3-2(t+1)^2+4H_{t},\
G^{HK}= \alpha-H_{T},\\&
Q_{t}^{HK}+\bar Q_{t}^{HK}= 4+\dot{K_{t}}+2K_{t}+\delta_{2t} H_{t},\
S_{t}^{HK}+\bar S_{t}^{HK}= 2+K_{t},\\&
R_{t}^{HK}+\bar R_{t}^{HK}= 1-2(t+1)^2+4H_{t},\\&
G^{HK}+\bar G^{HK}=-1-K_{T}, \ \ \  \forall H,K\in \Phi,
\end{aligned}\end{equation*}
where $\delta_{1t}=\int_{\Theta} E_{t,\theta}^2\nu(d\theta)> 0, \delta_{2t}=\int_{\Theta} (E_{t,\theta}+\bar E_{t,\theta})^2\nu(d\theta)>0$. In particular, if we define $H_{0 t}=\frac{1}{2}(t+1)^2, K_{0 t}=\frac{1}{1+(T-t)}-2$, then
\begin{equation*}
\begin{aligned}
&Q_{t}^{H_0K_0}= (t+1)+2(t+1)^2+\frac{\delta_1(t+1)^2}{2},
\ S_{t}^{H_0K_0}= \frac{1}{2}(t+1)^2,\\&
R_{t}^{H_0K_0}= (t+1)^3,\ G^{H_0K_0}= \alpha-\frac{1}{2}(T+1)^2,\\&
Q_{t}^{H_0K_0}+\bar Q_{t}^{H_0K_0}=  \frac{1}{(1+T-t)^2}+\frac{2}{1+T-t}+\frac{\delta_2(t+1)^2}{2},\\&
S_{t}^{H_0K_0}+\bar S_{t}^{H_0K_0}=  \frac{1}{1+(T-t)},\\&
R_{t}^{H_0K_0}+\bar R_{t}^{H_0K_0}= 1,\
G^{H_0K_0}+\bar G^{H_0K_0}= 0.
\end{aligned}
\end{equation*}
It is easy to see that Assumption (S) holds true for $J^{H_0K_0}[u]$. Then  Theorem \ref{main results Th} and Lemma \ref{H system Indefinite} imply that for any $H,K\in \Phi$, the   MF-FBSDEJ
\begin{equation*}
\left\{\begin{aligned}&dX_{t}^{HK}=\left(2X_{t}^{HK}-\mathbb E[X_{t}^{HK}]+u_{t}^{HK}\right)dt +2u_{t}^{HK}dW_{t}\\&\qquad\qquad+\int_{\Theta}( E_{t,\theta}X_{t-}^{HK}+\bar E_{t,\theta}\mathbb E[X_{t-}^{HK}])\widetilde N(dt,d\theta),
\\&dY_{t}^{HK}=-\left(2Y_{t}^{HK}-\mathbb E[Y_{t}^{HK}] +\int_{\Theta} (E_{t,\theta}r_{t,\theta}^{HK}+\bar E_{t,\theta}\mathbb E[r_{t,\theta}^{HK}])\nu(d\theta)\right.\\&\left.\qquad\qquad+Q_{t}^{HK}X_{t}^{HK}+\bar Q_{t}^{HK}\mathbb E[X_{t}^{HK}]+S_{t}^{HK}u_{t}^{HK}+\bar S_{t}^{HK}\mathbb E[u_{t}^{HK}]\right)dt\\&\qquad\qquad+Z_{t}^{HK}dW_{t}+\int_{\Theta} r_{t,\theta}^{HK}\widetilde N(dt,d\theta), \\&X_{0}^{HK}=x,\ \ \ \  Y^{HK}_{T}=G^{HK}X^{HK}_{T}+\bar G^{HK}\mathbb E[X^{HK}_{T}],\\
&S_{t}^{HK^\top} X_{t-}^{HK}+\bar S_{t}^{HK^\top}\mathbb E[X_{t-}^{HK}]+Y_{t-}^{HK}+2Z_{t}^{HK}+R_{t}^{HK}u_{t}^{HK}+\bar R_{t}^{HK}\mathbb E[u_{t}^{HK}]=0
\end{aligned}
\right.
\end{equation*}
 has a unique solution. Further, Theorem \ref{main results Th}  implies that Riccati equations
\begin{equation*}
\left\{\begin{aligned}&\dot{P_{t}}+4P_{t}+\delta_{1t} P_{t}-\frac{P_{t}^2}{(t+1)^3-2(t+1)^2+4P_{t}}=0,\\&P_{T}=\alpha
\end{aligned}\right.
\end{equation*}
and\begin{equation*}
\left\{\begin{aligned}&\dot{\Pi_{t}}+2\Pi_{t}+\delta_{2t}P_{t}+4-\frac{(\Pi_{t}+2)^2}{1-2(t+1)^2+4P_{t}}=0,\\ &\Pi_{T}=-1
\end{aligned}
\right.
\end{equation*}
admit unique  solutions $P, \Pi$, respectively. And  the
optimal pair $(X^*, u^*)$   satisfies
\begin{equation*}
\left\{\begin{aligned}
 u_{t}^*=&-\frac{P_{t}}{(t+1)^3-2(t+1)^2+4P_{t}}(X_{t}^*-\mathbb E[X_{t}^*])-\frac{\Pi_{t}}{1-2(t+1)^2+4P_{t}}\mathbb E[X_{t}^*],\\
dX_{t}^*=&\ \Big(2X_{t}^*-\mathbb E[X_{t}^*]+u_{t}^*\Big)dt+2u_{t}^*dW_{t}+\int_{\Theta}( E_{t,\theta}X_{t-}^*+\bar E_{t,\theta}\mathbb E[X_{t-}^*])\widetilde N(dt,d\theta),\\
X_{0}^*=&\ x.
\end{aligned}
\right.
\end{equation*}

We remark that the well-posedness of MF-FBSDE, i.e., the jump diffusion item in MF-FBSDEJ disappears, has been well studied (see \cite{Bensoussan2015Well,Carmona2012Mean,carmona2015forward,carmona2013probabilistic}). In detail, Bensoussan \cite{Bensoussan2015Well} derived the existence and uniqueness of solution of  MF-FBSDE  under a monotonicity condition.   Carmona and Delarue \cite{Carmona2012Mean}  obtained   the solvability of MF-FBSDE   by a compactness
argument and the Schauder fixed point theorem under a bound   condition. Carmona and Delarue \cite{carmona2015forward} took advantage
of the convexity of the Hamiltonian to apply the continuation method, and proved the existence and uniqueness of solution of MF-FBSDE. Carmona and Delarue \cite{carmona2013probabilistic}  derived the solvability results  by using an approximation procedure under some convexity condition. Moreover, Li and Min \cite{2020Fully} investigated the existence and uniqueness of solution  
to MF-FBSDEJ  under a monotonicity condition, which extended the results in previous literature. Different from the works above, our equivalent method provides an alternative way to solve MF-FBSDEJ. Specially, the original stochastic Hamiltonian system is of form
\begin{equation}\label{Application H1}
\left\{\begin{aligned}&dX_{t}=\left(2X_{t}-\mathbb E[X_{t}]+u_{t}\right)dt +2u_{t}dW_{t}+\int_{\Theta} (E_{t,\theta}X_{t-}+\bar E_{t,\theta}\mathbb E[X_{t-}])\widetilde N(dt,d\theta),\\&dY_{t}=-\left(2Y_{t}-\mathbb E[Y_{t}]+\int_{\Theta}( E_{t,\theta}r_{t,\theta}+\bar E_{t,\theta}\mathbb E[r_{t,\theta}])\nu(d\theta)\right)dt\\&\qquad\quad+Z_{t}dW_{t}+\int_{\Theta} r_{t,\theta}\widetilde N(dt,d\theta), \\&X_{0}=x,\ \ \ \  Y_{T}=\alpha X_{T}-(\alpha-1)\mathbb E[X_{T}],\\
&Y_{t}+2Z_{t}+R_{t}u_{t}+\bar R_{t}\mathbb E[u_{t}]=0.
\end{aligned}
\right.
\end{equation}
Since $R_{t}=(t+1)^3-2(t+1)^2,  R_{t}+\bar R_{t}=1-2(t+1)^2$, we can not derive an expression of the optimal control process $u$ from the last equation in \eqref{Application H1}. The  monotonicity condition in \cite{Bensoussan2015Well,2020Fully} and the bounded condition in \cite{Carmona2012Mean} fail. Moreover, Carmona and Delarue \cite{carmona2015forward,carmona2013probabilistic} assumed that cost functional satisfies some convex condition, which is not true in our setting, thus the methods in \cite{carmona2013probabilistic} and \cite{carmona2015forward} fail.
We emphasize that our method is also effective in proving the solvability of MF-FBSDEJ with a slightly general and complicated form. The following example provides   a better understanding on this issue.

{\it\textbf{Example 5.3}}:
Consider an  MF-FBSDEJ
\begin{equation}\label{Example 3}
\left\{\begin{aligned}dX_{t}=&\left(2X_{t}+\mathbb E[X_{t}]+Y_{t}+Z_{t}\right)dt+(Y_{t}+Z_{t})dW_{t}+\int_{\Theta}\bar E_{t,\theta}\mathbb E[X_{t-}]\widetilde N(dt,d\theta),\\dY_{t}=&-\left\{2Y_{t}+\mathbb E[Y_{t}]-X_{t}+\mathbb E[X_{t}]+\int_{\Theta} \bar E_{t,\theta}\mathbb E[r_{t,\theta}]\nu(d\theta)\right\}dt\\&+Z_{t}dW_{t}+\int_{\Theta}r_{t,\theta}\widetilde N(dt,d\theta),\\X(0)=&x,\ \ \ \  Y_{T}=2X_{T}-\mathbb E[X_{T}],
\end{aligned}
\right.
\end{equation}
where $X, Y, Z, r$ are 1-dimensional stochastic processes.
We claim that \eqref{Example 3} does not satisfy the conditions in \cite{Bensoussan2015Well,Carmona2012Mean,2020Fully}. Indeed, we have
\begin{equation*}
\begin{aligned}
&\mathbb E\left[(X_{1t}-X_{2t})\left(2Y_{2t}+\mathbb E[Y_{2t}]-X_{2t}+\mathbb E[X_{2t}]+\int_{\Theta}\bar E_{t,\theta}\mathbb E[r_{2t,\theta}]\nu(d\theta)\right.\right.\\&-\left.\left.\left(2Y_{1t}+\mathbb E[Y_{1t}]-X_{1t}+\mathbb E[X_{1t}]+\int_{\Theta}\bar E_{t,\theta}\mathbb E[r_{1t,\theta}]\nu(d\theta)\right)\right)\right]\\&+\mathbb E\bigg[(Y_{1t}-Y_{2t})\Big(2X_{1t}+\mathbb E[X_{1t}]+Y_{1t}+Z_{1t}-\left(2X_{2t}+\mathbb E[X_{2t}]+Y_{2t}+Z_{2t}\right)\Big)\bigg]\\&+\mathbb E\left[(Z_{1t}-Z_{2t}) \Big(Y_{1t}+Z_{1t}-\left(Y_{2t}+Z_{2t}\right)\Big)\right]\\&+\mathbb E\left[\int_{\Theta} (r_{1t,\theta}-r_{2t,\theta})\bar E_{t,\theta}\Big(\mathbb E[X_{1t-}]-\mathbb E[X_{2t-}]\Big)\nu(d\theta)\right]\\=&\ \mathbb E[(X_{1t}-X_{2t})^2]-\mathbb E[X_{1t}-X_{2t}]^2+\mathbb E[(Y_{1t}-Y_{2t}+Z_{1t}-Z_{2t})^2]\\ \geq &\ 0,\\
&\mathbb E\Big[(X_{1T}-X_{2 T})\Big(2X_{1 T}-\mathbb E[X_{1T}]-(2X_{2T}-\mathbb E[X_{2 T}])\Big)\Big]\\\geq &\ \mathbb E[(X_{1 T}-X_{2 T})^2].
\end{aligned}
\end{equation*}

It implies that  the monotonicity condition in \cite{Bensoussan2015Well,2020Fully} fails. \eqref{Example 3} is a linear MF-FBSDEJ, which does not satisfy the  bounded condition in \cite{Carmona2012Mean}.

We now prove the existence and uniqueness of solution of MF-FBSDEJ \eqref{Example 3} with the help of our results. Consider an  MF-LQJ problem with a 1-dimensional  state equation
\begin{equation*}\label{Example 3state}
\left\{\begin{aligned}dX_{t}=&\ \left(2X_{t}+\mathbb E[X_{t}]+u_{t}\right)dt+u_{t}dW_{t}+\int_{\Theta}\bar E_{t,\theta}\mathbb E[X_{t-}]\widetilde N(dt,d\theta),\\X_{0}=&\ x.
\end{aligned}
\right.
\end{equation*}
An  admissible control $u$ is  a predictable process such that $u\in L_{\mathbb F}^2(0,T;\mathbb R)$.
Introduce a cost functional
\begin{equation}\label{cost functional example3}
\begin{aligned}
J[u]= &\frac{1}{2}\mathbb E\left[2X_{T}^2-\mathbb E[X_{T}]^2+\int_0^T\left(-u_{t}^2-X_{t}^2+\mathbb E[X_{t}]^2\right)dt\right].
\end{aligned}
\end{equation}

The corresponding Hamiltonian system is
\begin{equation}\label{Application H3}
\left\{\begin{aligned}&dX_{t}=\left(2X_{t}+\mathbb E[X_{t}]+u_{t}\right)dt +u_{t}dW_{t}+\int_{\Theta}\bar E_{t,\theta}\mathbb E[X_{t-}]\widetilde N(dt,d\theta),\\&dY_{t}=-\left(2Y_{t}+\mathbb E[Y_{t}]-X_{t}+\mathbb E[X_{t}]+\int_{\Theta}\bar E_{t,\theta}\mathbb E[r_{t,\theta}]\nu(d\theta)\right)dt\\&\qquad\quad+Z_{t}dW_{t}+\int_{\Theta}r_{t,\theta}\widetilde N(dt,d\theta), \\&X_{0}=x,\ \ \ \  Y_{T}=2X_{T}-\mathbb E[X_{T}],\\
&-u_{t}+Y_{t-}+Z_{t}=0.
\end{aligned}
\right.
\end{equation}
Different from \eqref{Application H1},  we can  derive an explicit expression of the optimal control process $u$ from the last equation in \eqref{Application H3}.  In fact, it is easy to see that stochastic Hamiltonian system  \eqref{Application H3} is exactly  MF-FBSDEJ \eqref{Example 3} with $u_{t}=Y_{t}+Z_{t}$. Note that the cost functional does not satisfy the  convex conditions in \cite{carmona2015forward,carmona2013probabilistic}.

For the above MF-LQJ problem, it is clear  that Assumption (S) does  not hold. Now we introduce an equivalent cost functional $J^{H_0K_0}[u]$ satisfying Assumption (S). Recalling \eqref{cost HK}, we have
\begin{equation*}
\begin{aligned}
&Q_{t}^{HK}= -1+\dot{H_{t}}+4H_{t},\
S_{t}^{HK}= H_{t},\\&
R_{t}^{HK}= -1+H_{t},\
G^{HK}= 2-H_{T},\\&
Q_{t}^{HK}+\bar Q_{t}^{HK}= \dot{K_{t}}+6K_{t}+\delta_{t} H_{t},\
S_{t}^{HK}+\bar S_{t}^{HK}= K_{t},\\&
R_{t}^{HK}+\bar R_{t}^{HK}= -1+H_{t},\
G^{HK}+\bar G^{HK}= 1-K_{T}, \ \forall H,K\in \Phi,
\end{aligned}
\end{equation*}
where $\delta_{t}=\int_{\Theta}\bar E_{t,\theta}^2\nu(d\theta)$.
In particular, if we define $H_{0 t}=2, K_{0 t}=1$, then
\begin{equation*}
\begin{aligned}
&Q_{t}^{H_0K_0}= 7,\
S_{t}^{H_0K_0}= 2,\
R_{t}^{H_0K_0}= 1,\
G^{H_0K_0}= 0,\\&
Q_{t}^{H_0K_0}+\bar Q_{t}^{H_0K_0}= 6+2\delta_{t},\
S_{t}^{H_0K_0}+\bar S_{t}^{H_0K_0}=1,\\&
R_{t}^{H_0K_0}+\bar R_{t}^{H_0K_0}= 1,\
G^{H_0K_0}+\bar G^{H_0K_0}= 0.
\end{aligned}
\end{equation*}

With the data, Assumption (S) holds true for $J^{H_0K_0}[u]$. Then it follows  from Theorem \ref{main results Th} that MF-FBSDEJ \eqref{Application H3} admits a unique solution $(X, u, Y, Z,r)$, and $(X, Y, Z,r)$ is exactly the solution of \eqref{Example 3}.

{\it\textbf{Example 5.4}}:
Consider a financial market consisting of a bond and
a stock,  in which two assets are trading continuously within the time horizon $[0, T]$.  The dynamics of the bond price process $S_{1 t}$ is
governed by
\begin{equation*}
\left\{ \begin{aligned}
dS_{1t}=&\ r_{t}S_{1t}dt ,\\
S_{10}=&\ s_1,
\end{aligned}\right.
\end{equation*}
where $r_{t}$ is the interest rate of the bond. The dynamics of the stock price process $S_{2t}$ is
governed by
\begin{equation*}
\left\{ \begin{aligned}
dS_{2t}=&\ \mu_{t}S_{2 t}dt+\sigma_{t}S_{2t}dW_{t},\\
S_{20}=&\ s_2,
\end{aligned}\right.
\end{equation*}
where $\mu_{t}$ and $\sigma_{t}$ are the appreciation rate   and the volatility coefficient of the stock, respectively. For simplicity, we assume that the coefficients $\mu_{t}\geq r_{t}> 0$, $\sigma_{t}$ and $\frac{1}{\sigma_{t}}$ are bounded and deterministic functions.

We assume that the trading of shares takes place continuously in a self-financing fashion and there are no transaction costs. We denote by $N_{t}$ the asset of an investor and by $u_{t}$   the amount allocated in the stock share at time $t$. Clearly, the amount invested in the risk-free asset is $N_{t}-u_{t}$. Without liability, the
asset of the investor  $N_{t}$, evolves as
\begin{equation}\label{wealth}
\left\{ \begin{aligned}
dN_{t}=&\ [r_{t}N_{t}+(\mu_{t}-r_{t})u_{t}]dt+\sigma_{t}u_{t}dW_{t},\\
N_{0}=&\ n_0.
\end{aligned}\right.
\end{equation}
The investor's accumulative liability at time $t$  is denoted  by $L_{t}$. Chiu and Li \cite{chiu2006asset}, Wei and Wang \cite{wei2017time} described the liability process by  a geometric Brownian motion. In fact, it is possible that the control strategy and the mean of asset of the investor can influence the liability process, due to
the complexity of the financial market and the risk aversion
behavior of the investor. Such an example can be found
in Wang et al. \cite{Wang2017An}, where the liability process depends
on a control strategy (for example, capital injection or
withdrawal) of the firm.  Along this line, we proceed to improve the liability process here. Suppose that the dynamics of  $L_{t}$ satisfies
\begin{equation}\label{liability}
\left\{ \begin{aligned}
dL_{t}=&\ \big(a_{t}L_{t}+c_{t}\mathbb E[N_{t}]\big)dt+b_{t}L_{t}dW_{t} ,\\
L_{0}=&\ l_0,
\end{aligned}\right.
\end{equation}
where $a_{t}$ is the appreciation rate of  the liability  and $b_{t}$ is the corresponding volatility which satisfies the non-degeneracy condition. $a_{t}, b_{t}, c_{t}$ are deterministic continuous functions on $[0, T]$. Taking the liability into consideration, the SDE for the net wealth of the investor at time $t$, denote by $I_{t}$,  is obtained by subtracting \eqref{liability} from \eqref{wealth},
\begin{equation}\label{finance state}
\left\{ \begin{aligned}
dL_{t}=&\ \Big(a_{t}L_{t}+c_{t}\mathbb E[I_{t}]+c_{t}\mathbb E[L_{t}]\Big)dt+b_{t}L_{t}dW_{t},\\
L_0=&\ l_0,\\
dI_{t}=&\ \Big(r_{t}I_{t}+\big(r_{t}-a_{t}\big)L_{t}+\big(\mu_{t}-r_{t}\big)u_{t}-c_{t}\mathbb E[I_{t}]-c_{t}\mathbb E[L_{t}]\Big)dt+\big(\sigma_{t} u_{t}-b_{t}L_{t}\big)dW_{t},\\
I_{0}=&\ n_0-l_0.
\end{aligned}\right.
\end{equation}
\begin{definition}
An $\mathbb R$-valued portfolio strategy $u $ is called admissible,  if $u $ is $\mathbb F$-adapted and $\mathbb E\left[\int_0^Tu_{t}^2dt\right]<\infty$. The set of all admissible portfolio strategies is denoted by  $\mathcal U_{ad}$. 
\end{definition}

For any $u\in\mathcal  U_{ad}$, \eqref{finance state} admits a unique solution $(L,I)\in S_{\mathbb F}^2(0,T;\mathbb R^2)$. We introduce a performance functional of the investor, which is in the form of
$$
J[u]=\mathbb E\left\{\int_0^T\Big[L_{t}^2+(I_{t}-\mathbb E[I_{t}])^2\Big]dt+(I_{T}-\mathbb E[I_{T}])^2\right\}.$$
Now we pose   an asset-liability management problem as follows.

{\it\textbf {Problem AL}}: Find a portfolio strategy $u^*\in \mathcal U_{ad}$  such that
\begin{equation}\label{finance  inf}
J[u^*]=\inf_{u\in \mathcal U_{ad}}J[u].
\end{equation}
 The problem implies that the investor aims to minimize the  risk of the net wealth  and  the liability over the whole time horizon, simultaneously.

It is easy to see that Problem AL is a special case of Problem MF. Denoting $X_{t}=(L_{t}, I_{t})^\top$, and $x=(l_0, n_0-l_0)^\top$, we have
\begin{equation*}
\begin{aligned}
&A_{t}=\left(\begin{array}{cc}
a_{t}&0\\
r_{t}-a_{t}&r_{t}\end{array}\right),\  \bar A_{t}=\left(\begin{array}{cc}
c_{t}&c_{t}\\
-c_{t}&-c_{t}\end{array}\right),B_{t}=\left(\begin{array}{c}0\\
\mu_{t}-r_{t}\end{array}\right),
C=\left(\begin{array}{cc}
b_{t}&0\\
-b_{t}&0\end{array}\right),\\& D=\left(\begin{array}{c}0\\
\sigma_{t}\end{array}\right),\ Q=\left(\begin{array}{cc}
1&0\\
0&1\end{array}\right),\ \bar Q=\left(\begin{array}{cc}
0&0\\
0&-1\end{array}\right), G=\left(\begin{array}{cc}
0&0\\
0&1\end{array}\right),
\bar G=\left(\begin{array}{cc}
0&0\\
0&-1\end{array}\right),\\&\bar B=0_{2\times1},\ \bar C=0_{2\times2},\ \bar D=0_{2\times1},\ R=\bar R=0,\ S=\bar S=0_{2\times1}.
\end{aligned}\end{equation*}
Clearly, Assumption (S) does not hold. Define $H_{t}=\left(\begin{array}{cc}
0&0\\
0&\lambda_{t}\end{array}\right), K_{t}=0_{2\times2}$, where $\lambda_{t}$ is the solution of
\begin{equation*}
\left\{ \begin{aligned}
\dot\lambda_{t}=&\ \left[\left(\frac{\mu_{t}-r_{t}}{\sigma_{t}}\right)^2-2r_{t}\right]\lambda_{t}+\left[r_{t}-a_{t}+\frac{b_{t}(\mu_{t}-r_{t})}{\sigma_{t}}\right]^2\lambda_{t}^2,\\
\lambda_{T}=&\ \frac{1}{2}.
\end{aligned}\right.
\end{equation*}
Recalling \eqref{cost HK}, we have

\begin{equation*}
\begin{aligned}
&Q^{HK}_{t}=\left(\begin{array}{cc}
1+b_{t}^2\lambda_{t}&(r_{t}-a_{t})\lambda_{t}\\
(r_{t}-a_{t})\lambda_{t}&1+2r_{t}\lambda_{t}+\dot \lambda_{t}\end{array}\right),S^{HK}_{t}=\left(\begin{array}{c}
-b_{t}\sigma_{t}\lambda_{t}\\
(\mu_{t}-r_{t})\lambda_{t}\end{array}\right),\\&\
R^{HK}_{t}=\lambda_{t}\sigma_{t}^2,G^{HK}=\left(\begin{array}{cc}
0&0\\
0&\frac{1}{2}\end{array}\right),\ Q^{HK}_{t}+\bar Q^{HK}_{t}=\left(\begin{array}{cc}
1+b_{t}^2\lambda_{t}&0\\
0&0\end{array}\right), \\&S^{HK}_{t}+\bar S^{HK}_{t}=\left(\begin{array}{c}
-b_{t}\sigma_{t}\lambda_{t}\\
0\end{array}\right),R^{HK}_{t}+\bar R^{HK}_{t}=\lambda_{t}\sigma_{t}^2,\ G^{HK}+\bar G^{HK}=0.
\end{aligned}
\end{equation*}
It is easy to see that Assumption  (S) holds true for   $J^{HK}[u]$. According to Theorem \ref{main results Th}, we know that  the following Riccati equations have  unique solutions
\begin{equation}\label{P riccati finance}
\left\{\begin{aligned}&\dot{P_{t}}+P_{t}A_{t}+A_{t}^\top P_{t}+C_{t}^\top P_{t}C_{t}+Q_{t}-(P_{t}B_{t}+C_{t}^\top P_{t}D_{t})\Sigma_{t}^{-1}(B_{t}^\top P_{t}+D_{t}^\top P_{t}C_{t})=0,\\&P_{T}=G,
\end{aligned}\right.
\end{equation}
\begin{equation}\label{Pi riccati finance}
\left\{\begin{aligned}&\dot{\Pi_{t}}+\Pi_{t}(A_{t}+\bar A_{t})+(A_{t}+\bar A_{t})^\top \Pi_{t}+C_{t}^\top P_{t}C_{t}+(Q_{t}+\bar Q_{t})\\&\ \ -[\Pi_{t}B_{t}+C_{t}^\top P_{t}D_{t}]\Sigma_t^{-1}[B_{t}^\top \Pi_{t}+D_{t}^\top P_{t}C_{t}]=0,\\ &\Pi_{T}=G+\bar G,
\end{aligned}
\right.
\end{equation}
where
\begin{equation*}\Sigma_{t}=D_{t}^\top P_{t}D_{t}.\end{equation*}
An optimal portfolio strategy $u^*$ is given by
\begin{equation}\label{financecontrol}
\begin{aligned}u_{t}^*=&-\Sigma_{t}^{-1}(B_{t}^\top P_{t}+D_{t}^\top P_{t}C_{t})(X_{t}-\mathbb E[X_{t}])-\Sigma_{t}^{-1}[B_{t}^\top \Pi_{t}+D_{t}^\top P_{t}C_{t}]\mathbb E[X_{t}],
\end{aligned}
\end{equation}
where
\begin{equation*}
\left\{\begin{aligned}dX_{t}=&\Big\{\left[A_{t}-B_{t}\Sigma_{t}^{-1}(B_{t}^\top P_{t}+D_{t}^\top P_{t}C_{t})\right](X_{t}-\mathbb E[X_{t}])\\&+\left[A_{t}+\bar A_{t}-B_{t}\Sigma_{t}^{-1}(B_{t}^\top \Pi_{t}+D_{t}^\top P_{t}C_{t})\right]\mathbb E[X_{t}]\Big\}dt\\&+\Big\{\left[C_{t}-D_{t}\Sigma_{t}^{-1}(B_{t}^\top P_{t}+D_{t}^\top P_{t}C_{t})\right](X_{t}-\mathbb E[X_{t}])\\&+\left[C_{t}-D_{t}\Sigma_{t}^{-1}(B_{t}^\top \Pi_{t}+(D_{t}^\top P_{t}C_{t})\right]\mathbb E[X_{t}]\Big\}dW_{t}, \\
X_{0}=&\ x.
\end{aligned}
\right.
\end{equation*}

Note that it is hard to give  a more explicit expression of $u^*_{t}$ due to the complexity of  Riccati equations \eqref{P riccati finance} and \eqref{Pi riccati finance}. We will use numerical simulation to illustrate the   optimal investment strategy $u^*_{t}$ and to analyse the relationship between $u^*_{t}$  and some common financial parameters in our model. Here we only  analyse the relationship between $u^*_{t}$ and the risk-free rate $r_{t}$, as well as the appreciation rate of the liability $a_{t}$, respectively. In the following discussion, we fix $T = 1, n_0 =1, l_0=0.5, \mu_{t}=0.3, \sigma_{t}=0.5, b_{t}=0.3, c_{t}=0.5$.

\begin{figure}[htbp]
\centering
\includegraphics[height=3.5in,width=7in]{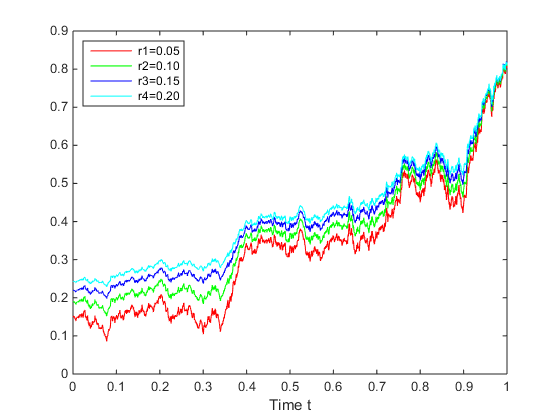}
\caption{The relationship between $r$ and $u^*_{t}$.}
\end{figure}

\begin{figure}[htbp]
\centering
\includegraphics[height=3.5in,width=7in]{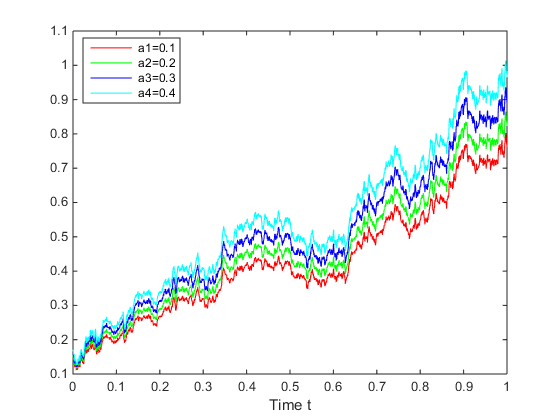}
\caption{The relationship between $a$ and $u^*_{t}$.}
\end{figure}

{\it\textbf{A.}}: The relationship between $r$ and $u^*_{t}$.
\\Let $a_{t}=0.2$, $[r_{1 t}, r_{2 t}, r_{3 t}, r_{4 t}]=[0.05, 0.1, 0.15, 0.2]$. We plot Fig. 5.1 which illustrates the relationship between the optimal investment strategy $u^*_{t}$ and the risk-free rate $r_{t}$.
From Fig. 5.1, we   find that the higher the risk-free rate is, the higher the amount of investor's wealth allocated in the stock share is. It is   reasonable that the net wealth expectation $\mathbb E[I_{t}]$ increases    with the increase of the risk-free rate, which leads to the increase of the liability $L_{t}$. Thus, the investor has to invest more in stock share to avoid the  risk of the net wealth.

{\it\textbf{B}}: The relationship between $a$ and $u^*_{t}$.\\
Let $r_{t}=0.05$, $[a_{1 t}, a_{2 t}, a_{3 t}, a_{4 t}]=[0.1, 0.2, 0.3, 0.4]$. We plot Fig. 5.2 which  illustrates the relationship between the optimal investment strategy $u^*_{t}$ and the appreciation rate of the liability $a_{t}$.
 Fig. 5.2 shows that the higher the appreciation rate of the liability is, the higher the amount of investor's wealth allocated in the stock share is.  It is  reasonable that the liability   increases   with the increase of the  appreciation rate of liability.  Consequently,  the investor has to invest more in stock share to avoid the  risk of the net wealth.

\section{Concluding remarks}
We use an equivalent cost functional method to deal with an indefinite MF-LQJ problem. Relying on this method, we transform an indefinite MF-LQJ problem to a definite MF-LQJ problem.    Compared with existing literature, we further consider the unique  solvabilities of   Riccati equations arising in indefinite MF-LQJ problems, which have not  been  studied before. Moreover,  the method provides an alternative approach to solve MF-FBSDEJ. Several illustrative
examples show that our method is more effective than existing
methods in proving the existence and uniqueness of solution to
MF-FBSDEJ.
Our results are obtained with the framework of MF-LQJ problem. We will further investigate some results for other stochastic system. For example,  LQ optimal control  with delay has not been completely solved yet. We will try to extend this method to the framework with delay.

\end{document}